\documentclass[12pt]{article}

\usepackage{amsmath,amssymb}

\newtheorem{theorem}{Theorem}

\newtheorem{corollary}[theorem]{Corollary}

\newtheorem{lemma}[theorem]{Lemma}

\newtheorem{proposition}[theorem]{Proposition}
\newtheorem{remark}[theorem]{Remark}

\begin{document}

\title{\textbf{
Tubes estimates for diffusion processes\\under a local H\"ormander condition\\
of order one}}
\author{ \textsc{Vlad Bally}\thanks{%
Laboratoire d'Analyse et de Math\'ematiques Appliqu\'ees, UMR 8050,
Universit\'e Paris-Est Marne-la-Vall\'ee, 5 Bld Descartes, Champs-sur-Marne,
77454 Marne-la-Vall\'ee Cedex 2, France. Email: \texttt{bally@univ-mlv.fr} }%
\smallskip \\
%EndAName
\textsc{Lucia Caramellino}\thanks{%
Dipartimento di Matematica, Universit\`a di Roma - Tor Vergata, Via della
Ricerca Scientifica 1, I-00133 Roma, Italy. Email: \texttt{%
caramell@mat.uniroma2.it}}\smallskip\\
}
%\date{}
\maketitle

\parindent 0pt

{\textbf{Abstract.}}
We consider a diffusion process $X_{t}$ and a skeleton curve $x_{t}(\phi )$
and we give a lower bound for $P(\sup_{t\leq T}d(X_{t},x_{t}(\phi ))\leq R).$
This result is obtained under the hypothesis that the strong H\"{o}rmander
condition of order one (which involves the diffusion vector fields and the
first Lie brackets) holds in every point $x_{t}(\phi ),0\leq t\leq T.$ Here $%
d$ is a distance which reflects the non isotropic behavior of the diffusion
process which moves with speed $\sqrt{t}$ in the directions of the diffusion
vector fields but with speed $t$ in the directions of the first order Lie
brackets. We prove that $d$ is locally equivalent with the standard control
metric $d_{c}$ and that our estimates hold for $d_{c}$ as well.

\bigskip


\medskip


\clearpage

\tableofcontents

\section{Introduction}

We consider the diffusion process solution of $dX_{t}=\sum_{j=1}^{d}\sigma
_{j}(t,X_{t})\circ dW_{t}^{j}+b(t,X_{t})dt$ where the coefficients $\sigma
_{j},b$ are three times differentiable and verify the strong H\"{o}rmander
condition on order one (involving $\sigma _{j}$ and the first order Lie brackets  $[\sigma _{i},\sigma_{j}]$) locally around a skeleton path $dx_{t}(\phi )=\sum_{j=1}^{d}\sigma
_{j}(t,x_{t}(\phi ))\phi _{t}^{j}dt+b(t,x_{t}(\phi ))dt.$ The aim of this
paper is to give a lower bound for the probability that $X_{t}$ remains in a
tube around $x_{t}(\phi )$ for $t\leq T.$

This problem has already been addressed in the literature. The first result
was given by Stroock and Varadhan in their celebrated paper \cite{[SV]}. They
obtain a lower bound for $P(\sup_{t\leq T}\left\Vert X_{t}-x_{t}(\phi
)\right\Vert \leq R)$ and use it in order to prove the support theorem for
diffusion processes. Here $\left\Vert X_{t}-x_{t}(\phi )\right\Vert $ is the
Euclidian norm. Later, one has considered other norms which reflect the
degree of regularity of the trajectories of the diffusion process $X_{t}$:
Ben Arous and Gradinaru \cite{[BG]} and Ben Arous, Gradinaru and Ledoux \cite{[BGL]}
obtained similar results for the H\"{o}lder norm. And more recently Friz,
Lyons and Stroock \cite{[FLS]} use a norm related to the rough path theory. All
these results hold without any non degeneracy assumption.

Tubes estimates has also been considered in connection with the
Onsager-Machlup functional for diffusion processes. There is an abundant
literature on this subject: see e.g. \cite{[C1]}, \cite{[C2]}, \cite{[HT]}, \cite{[IW]}, \cite{[Z]}. In this case one
considers strong ellipticity conditions and the norm which describes the
tube is the Euclidian norm or some H\"{o}lder norm. Notice that these are
asymptotic results whether in our paper we give estimates which are non
asymptotic. Finally, in \cite{[BFM]} and \cite{[BM]} one obtains similar lower bounds for
general It\^{o} processes under an ellipticity assumption.

\smallskip

The specific point in our paper is that we use a distance which reflects the
non isotropic structure of the problem: the diffusion process $X_{t}$ moves
with speed $\sqrt{t}$ in the direction of the diffusion vector fields $%
\sigma _{j}$ and with speed $t=\sqrt{t}\times \sqrt{t}$ in the direction of $%
[\sigma _{i},\sigma _{j}].$ Let us be more precise. For $R>0$ and $x\in R^{n}$ we construct the matrix $A_{R}(t,x)$ with columns $\sqrt{R}\sigma _{i}(t,x),[%
\sqrt{R}\sigma _{j},\sqrt{R}\sigma _{p}](t,x),1\leq i,j,p\leq d.$ If the above
vectors span $R^{n}$ the matrix $A_{R}A_{R}^{\ast }(t,x)$ is invertible, so
we are able to define the norm
\begin{equation*}
\left\vert y\right\vert _{A_{R}(t,x)}^{2}=\left\langle (A_{R}A_{R}^{\ast
})^{-1}(t,x)y,y\right\rangle .
\end{equation*}%
Our main result is the following (see Theorem \ref{NOT3} for a precise
statement): we assume that the non-degeneracy condition holds along the
curve $x_{t}(\phi ),0\leq t\leq T$ and we prove
\begin{equation*}
P(\sup_{t\leq T}\left\vert X_{t}-x_{t}(\phi )\right\vert _{A_{R}(t,x_{t}(\phi
))}\leq 1)\geq \exp \Big(-C\Big(\frac{1}{R}+\int_{0}^{T}\left\vert \phi
_{t}\right\vert ^{2}dt\Big)\Big).
\end{equation*}%
Computations involving the above norms are generally not easy - so we give
some estimates which seem to be more explicit. In Proposition \ref{NOT1} we
prove that $\left\vert y\right\vert _{A_{R}(t,x)}$ describes (roughly
speaking) ellipsoids with semi-axes of length $\sqrt{R}$ in the directions
of $\sigma _{j}(t,x)$ and of length $R$ in the directions of $[\sigma
_{i},\sigma _{j}](t,x).$ Moreover we associate to the above norms the
following semi-distance: $d(x,y)<R$ if and only if $\left\vert y\right\vert
_{A_{R}(x)}<1.$ With this definition we have $\{\sup_{t\leq T}\left\vert
X_{t}-x_{t}(\phi )\right\vert _{A_{R}(t,x_{t}(\phi ))}\leq 1\}=\{\sup_{t\leq
T}d(x_{t}(\phi ),X_{t})\leq R\}.$ In Proposition \ref{NORM4} we prove that
the semi-distance $d$ is equivalent with the standard control metric $d_{c}$
(see (\ref{Not14}) for the definition) so the estimates of the tubes hold in
the control metric as well. In Proposition \ref{NOT2} we give local lower
and upper bounds for $d$ and $d_{c}$ in terms of some semi-distances which
describe in a more explicit way the ellipsoid structure we mentioned above.

The paper is organized as follows. In Section 2 we give the statements of
the main results. In Section 3 we consider a process $Z_{t}$ which is a
linear combination of $W_{t}^{j},j=1,...,d$ and of $%
\int_{0}^{t}W_{s}^{i}dW_{s}^{j},1\leq i,j\leq d.$ And we give a
decomposition of such a process - this decomposition represents the main
ingredient in our approach. Roughly speaking the idea is the following: we
consider a small interval of time $[0,\delta ]$ and we split it in $d$ subintervals $I_{i}=(t_{i-1},t_{i}]$ with $t_{i}=\frac{i}{d}\delta .$ We fix $i$
and for $t\in I_{i}$ we take conditional expectation with respect to $%
W_{t}^{j},j\neq i$ so all these processes appear as ``controls''. And the only
process which is at work is $W_{t}^{i}.$ Then the vector $%
(W_{t_{i}}^{i}-W_{t_{i-1}}^{i}),%
\int_{t_{i-1}}^{t_{i}}(W_{s}^{j}-W_{t_{i-1}}^{j})dW_{s}^{i},j\neq i$ is
Gaussian (with respect to the above mentioned conditional probability). And
we may choose the trajectories (controls) $(W_{s}^{j}-W_{t_{i-1}}^{j})_{s\in
I_{i}},j\neq i$ in such a way that the covariance matrix of the above
Gaussian vector is non degenerated (this is a support property proven in
Section 7). Then we are able to use estimates for non degenerated Gaussian
random variables. The process $Z_{t}$ appears as the principal part in the
development in stochastic series of order two of the diffusion process $%
X_{t}.$ In Section 4 we use the estimates for $Z_{t}$ in order to obtain
estimates for $X_{t}$ and so to finish the proof of the main theorem stated
in Section 2.

The fact that one may choose $(W_{s}^{j}-W_{t_{i-1}}^{j})_{s\in I_{i}},j\neq
i$ in an appropriate way is due to the support theorem for the Brownian
motion. But the quantitative property that we use employs in a crucial way
the estimates of the variance (with respect to the time) of the Brownian
motion obtained in \cite{[DY]}.

\bigskip

\textbf{Acknowledgments}. We are grateful to Arturo Kohatsu-Higa and to
Peter Friz for useful discussions on this topic.

\section{Notations and main results}

We consider the $n$ dimensional diffusion process
\begin{equation}
dX_{t}=\sum_{j=1}^{d}\sigma _{j}(t,X_{t})\circ dW_{t}^{j}+b(t,X_{t})dt
\label{Not1}
\end{equation}%
where $W=(W^{1},...,W^{d})$\ is a standard Brownian motion, $\circ
dW_{t}^{j} $ denotes the Stratonovich integral and $\sigma
_{j},b:R_{+}\times R^{n}\rightarrow R^{n}$ are three time differentiable in $%
x\in R^{n}$ and one time differentiable with respect to the time $t\in
R_{+}. $ We also assume that the derivatives with respect to the space $x\in
R^{n}$ are one time differentiable with respect to $t.$ And for $(t,x)\in
R_{+}\times R^{n}$ we denote by $n(t,x)$ a constant such that for every $%
s\in \lbrack (t-1)\vee 0,t+1],y\in B(x,1)$ and for every multi index $%
\alpha $ of length less or equal to three
\begin{equation}
\left\vert \partial _{x}^{\alpha }b(s,y)\right\vert +\left\vert \partial
_{t}\partial _{x}^{\alpha }b(s,y)\right\vert +\sum_{j=1}^{d}\left\vert
\partial _{x}^{\alpha }\sigma _{j}(s,y)\right\vert +\left\vert \partial
_{t}\partial _{x}^{\alpha }\sigma _{j}(s,y)\right\vert )\leq n(t,x).
\label{Not2}
\end{equation}%
Here, $\alpha =(\alpha _{1},...,\alpha _{k})\in \{1,...,n\}^{k}$ is a multi
index and $\left\vert \alpha \right\vert =k$ is the length of $\alpha$ and
$\partial _{x}^{\alpha }=\partial _{x_{\alpha _{1}}}...\partial _{x_{\alpha
_{k}}}.$ 
%Notice that we do not assume Lipschitz continuity and/or sublinear growth properties  for the coefficients, so it is not clear that the above $SDE$ has a unique solution.

\smallskip

In the following we assume that for external reasons one produces a
continuous adapted process $X$ which solves equation (\ref{Not1}) on the time interval $[0,T]$ and we give
estimates for this process.
More precisely, for $\phi \in L^{2}([0,T];R^{d})$, we assume there exists a solution of
\begin{equation}
dx_{t}(\phi )=\sum_{j=1}^{d}\sigma _{j}(t,x_{t}(\phi ))\phi
_{t}^{j}dt+b(t,x_{t}(\phi ))dt  \label{Not3}
\end{equation}%
and we want to estimate the probability that $X_{t}$ remains in a tube around
the deterministic curve $x_t=x_t(\phi)$.

\smallskip

We need some more notations. First, we use the following notation of
directional derivatives: for $f,g:R_{+}\times R^{n}\rightarrow R^{n}$ we
define $\partial _{g}f(t,x)=\sum_{i=1}^{n}g^{i}(t,x)\partial _{x_{i}}f(t,x)$
and we recall that the Lie bracket (with respect to the space variable $x)$
is defined as $[f,g](t,x)=\partial _{g}f(t,x)-\partial _{f}g(t,x).$
Moreover, let $M\in \mathcal{M}_{n\times m}$ \ be a matrix
(which generally may be not square) such that $MM^{\ast }$ is invertible ($%
M^{\ast }$ denotes the transposed matrix). We denote by $\lambda _{\ast
}(M)$ (respectively $\lambda ^{\ast }(M))$ the smaller (respectively the
larger) eigenvalue of $MM^{\ast }$ and we consider the norm on $R^{n}$
\begin{equation}
\left\vert y\right\vert _{M}=\sqrt{\left\langle (MM^{\ast
})^{-1}y,y\right\rangle }.  \label{Not4}
\end{equation}

We are concerned with the matrix $A(t,x)\in \mathcal{M}_{n\times m}$ with
columns $\sigma _{i}(t,x),[\sigma _{j},\sigma _{p}](t,x),1\leq i,j,p\leq
d,j\neq p.$ Here and all along the paper
\begin{equation*}
m=d^{2}.
\end{equation*}
We will write%
\begin{equation}
A(t,x)=(\sigma _{i}(x),[\sigma _{j},\sigma _{p}](t,x))_{i,j,p=1,...,d,j\neq
p}.  \label{Not9}
\end{equation}%
We denote by $\lambda (t,x)$ the lower eigenvalue of $A(t,x)$ that is
\begin{equation}
\lambda (t,x)=\inf_{\left\vert \xi \right\vert =1}\sum_{i=1}^{m}\left\langle
A_{i}(t,x),\xi \right\rangle ^{2},  \label{Not9'}
\end{equation}%
$A_{i}(t,x)$, $i=1,\ldots,m$, denoting the columns of $A(t,x)$. Moreover for $R>0$ we define
\begin{equation*}
A_{R}(t,x)=(\sqrt{R}\sigma _{i}(t,x),[\sqrt{R}\sigma _{j},\sqrt{R}\sigma
_{p}](t,x))_{i,j,p=1,...,d,j\neq p}.
\end{equation*}%
Consider now some $x\in R^{n},t\geq 0$ such that $(\sigma _{i}(t,x),[\sigma
_{j},\sigma _{p}](t,x))_{i,j,p=1,...,d,j\neq p}$ span $R^{n}.$ Then $%
A_{R}A_{R}^{\ast }(t,x)$ is invertible and we may define $\left\vert
y\right\vert _{A_{R}(t,x)}.$ We give some lower and upper bounds for $%
\left\vert y\right\vert _{A_{R}(t,x)}.$ We denote by $S(t,x)$ the space
spanned by $\sigma _{1}(t,x),...,\sigma _{d}(t,x)$ and by $S^{\bot }(t,x)$ the
orthogonal of $S(t,x).$ We also denote by $\Pi _{t,x}$ the projection on $%
S(t,x)$ and by $\Pi _{t,x}^{\bot }$ the projection on $S^{\bot }(t,x).$
Moreover we denote%
\begin{equation}
\lambda _{t,x}=\inf_{\xi \in S(t,x),\left\vert \xi \right\vert
=1}\sum_{i=1}^{d}\left\langle \sigma _{i}(t,x),\xi \right\rangle ^{2},\quad
\lambda _{t,x}^{\bot }=\inf_{\xi \in S^{\bot }(t,x),\left\vert \xi
\right\vert =1}\sum_{i<j}\left\langle [\sigma _{i},\sigma _{j}](t,x),\xi
\right\rangle ^{2}.  \label{Not8}
\end{equation}%
By the very definition $\lambda _{t,x}>0$ (which is different from $\lambda
(t,x))$ and under our hypothesis $\lambda _{t,x}^{\bot }>0$\ also. Then
Proposition \ref{NORM2} gives:

\begin{proposition}
\label{NOT1}If $R\leq \lambda _{t,x}/(4m\times n^{4}(t,x))$ then%
\begin{equation}
\frac{1}{4Rn^{2}(t,x)}\left\vert \Pi _{t,x}y\right\vert ^{2}+\frac{1}{%
4R^{2}n^{2}(t,x)}\left\vert \Pi _{t,x}^{\bot }y\right\vert ^{2}\leq
\left\vert y\right\vert _{A_{R}(t,x)}^{2}\leq \frac{4}{R\lambda _{t,x}}%
\left\vert \Pi _{t,x}y\right\vert ^{2}+\frac{4}{R^{2}\lambda _{t,x}^{\bot }}%
\left\vert \Pi _{t,x}^{\bot }y\right\vert ^{2}.  \label{Not10}
\end{equation}
\end{proposition}

For $\mu \geq 1$ and $0<h\leq 1$ we denote by $L(\mu ,h)$ the class of non
negative functions $f:R_{+}\rightarrow R_{+}$ which have the property%
\begin{equation*}
f(t)\leq \mu f(s)\quad \mbox{for}\quad \left\vert t-s\right\vert \leq h.
\end{equation*}%
We will make the following hypothesis: there exists some functions $%
n:[0,T]\rightarrow \lbrack 1,\infty )$ and $\lambda :[0,T]\rightarrow (0,1]$
such that for some $\mu \geq 1$ and $0<h\leq 1$ we have
\begin{equation}\label{Not6}
\begin{array}{lrl}
(H_1)\quad&n(t,x_{t}(\phi )) &\leq n_{t},\forall t\in \lbrack 0,T],\smallskip\\
(H_2)\quad&\lambda (t,x_{t}(\phi )) &\geq \lambda _{t}>0,\forall t\in
\lbrack 0,T],\smallskip\\
(H_3)\quad&n_{.},\lambda _{.} &\in L(\mu ,h).
\end{array}
\end{equation}

\begin{remark}
The hypothesis $(H_{2})$ implies that for each $t\in (0,T),$ the space $R^n$ is spanned by the vectors $(\sigma
_{i}(t,x_{t}),[\sigma _{j},\sigma _{p}](t,x_{t}))_{i,j,p=1,...,d,j<p}$, so the H\"{o}rmander condition holds along the curve $x_{t}(\phi ).$
\end{remark}

The main result in this paper is the following.

\begin{theorem}
\label{NOT3}Suppose that $(H_1)$, $(H_2)$ and $(H_3)$ hold and that $X_{0}=x_{0}(\phi ).$
Let $\rho \in (0,1).$ There exists a universal constant $C$ (depending on $d$
and $\rho $ only) such that for every $R\in (0,1)$ one has%
\begin{equation}
P(\sup_{t\leq T}\left\vert X_{t}-x_{t}(\phi )\right\vert
_{A_{R}(t,x_{t}(\phi ))}\leq 1)\geq \exp \Big(-C\mu ^{9}\Big(\frac{T}{h}+\int_{0}^{T}%
\frac{n_{t}^{6(1+d\rho )}}{\lambda _{t}^{1+2d\rho }}\Big(\frac{1}{R}+\left\vert
\phi _{t}\right\vert ^{2}\Big)dt\Big)\Big).  \label{Not7}
\end{equation}
\end{theorem}

\begin{remark}
Suppose that $X_{t}=W_{t}$ is just the Brownian motion and that $x_{t}=0$, so that $n_{t}=1$, $\lambda _{t}=1$, $\mu =1$ and $\phi _{t}=0$. Then $\left\vert
X_{t}-x_{t}\right\vert _{A_{R}(x_{t}(\phi ))}=R^{-1/2}W_{t}$ so we obtain $%
P(\sup_{t\leq T}\left\vert W_{t}\right\vert \leq \sqrt{R})\geq \exp (-CT/R)$
which is coherent with the standard estimate (see \cite{[IW]}).
\end{remark}

\begin{remark}
Since $\partial _{t}x_{t}(\phi )-b(t,x_{t}(\phi ))=\sigma (t,x_{t}(\phi
))\phi (t)$ we immediately obtain%
\begin{equation*}
\frac{1}{dn(t,x_{t}(\phi ))}\left\vert \partial _{t}x_{t}(\phi
)-b(t,x_{t}(\phi ))\right\vert \leq \left\vert \phi (t)\right\vert \leq
\frac{1}{\sqrt{\lambda _{t,x_{t}(\phi )}}}\left\vert \partial _{t}x_{t}(\phi
)-b(t,x_{t}(\phi ))\right\vert
\end{equation*}%
with $\lambda _{t,x_{t}(\phi )}$ given in (\ref{Not8}).
\end{remark}

We establish now the link between the norm $\left\vert z\right\vert
_{A_{R}(t,x)}$ and the control (Caratheodory) distance. We will use in a
crucial way the alternative characterizations given in \cite{[NSW]} for this
distance - and these results hold in the homogeneous case: the coefficients
of the equations do not depend on time: $\sigma _{j}(t,x)=\sigma
_{j}(x)$ and $b(t,x)=b(x).$ Consequently now on we have a matrix $A_{R}(x)$
instead of $A_{R}(t,x).$ We define the semi-distance $d:R^{n}\times
R^{n}\rightarrow R_{+}$ by $d(x,y)<\sqrt{R}$ if and only if $\left\vert
y\right\vert _{A_{R}(x)}<1$ (see page \pageref{Norm7} for the definition of a
semi-distance).

We also consider the standard control distance $d_{c}$ (Caratheodory
distance) associated to $\sigma _{1},...,\sigma _{d}$ in the following way.
Let $y_{t}(\phi )$ be the solution of the equation $dy_{t}(\phi
)=\sum_{j=1}^{d}\sigma _{j}(y_{t}(\phi ))\phi _{t}^{j}dt$ (notice that here $%
b=0).$ We denote $C(x,y)=\{\phi \in L^{2}(0,1):y_{0}(\phi )=x,y_{1}(\phi
)=y\}$ and we define%
\begin{equation}
d_{c}(x,y)=\inf \Big\{\Big(\int_{0}^{1}\left\vert \phi _{s}\right\vert
^{2}ds\Big)^{1/2}:\phi \in C(x,y)\Big\}.  \label{Not14}
\end{equation}%
In Section 8 Theorem \ref{NORM4} we prove that $d$ is \ locally equivalent
with $d_{c}$. Moreover we obtain the following bounds for them$.$ We define $%
\overline{d}(x,y)$ and $\underline{d}(x,y)$ as follows:

\begin{itemize}
\item
$\overline{d}(x,y)<\sqrt{R}$ if and only if
\begin{equation*}
\frac{4}{R\lambda _{x}}\left\vert \Pi _{x}(y-x)\right\vert ^{2}+\frac{4}{%
R^{2}\lambda _{x}^{\bot }}\left\vert \Pi _{x}^{\bot }(y-x)\right\vert ^{2}<1;
\end{equation*}%
\item
$\underline{d}(x,y)<\sqrt{R}$ if and
only if
\begin{equation*}
\frac{1}{4Rn_{x}^{2}}\left\vert \Pi _{x}(y-x)\right\vert ^{2}+\frac{1}{%
4R^{2}n_{x}^{2}}\left\vert \Pi _{x}^{\bot }(y-x)\right\vert ^{2}<1.
\end{equation*}%
\end{itemize}
Then as an immediate consequence (we give a detailed proof at the end of
Appendix 4) of Proposition \ref{NOT1} and Theorem \ref{NORM4} we obtain:

\begin{proposition}
\label{NOT2} Let $x,y\in R^{n}$ be such that%
\begin{equation}
\left\vert y-x\right\vert \leq \frac{\lambda _{x}\sqrt{\lambda _{\ast }(A(x))%
}}{(4m)n^{4}(x)}.  \label{Not11'}
\end{equation}%
Then
\begin{equation}
\underline{d}(x,y)\leq d(x,y)\leq \overline{d}(x,y).  \label{Not11}
\end{equation}%
Moreover for every compact set $K\subset R^{n}$ there exists some constants $%
C_{K},r_{K}$ such that for ever $x,y\in K$ which satisfy (\ref{Not11'}) and
such that $\overline{d}(x,y)\leq r_{K}$ one has
\begin{equation}
\frac{1}{C_{K}}\underline{d}(x,y)\leq d_{c}(x,y)\leq C_{K}\overline{d}(x,y).
\label{Not12}
\end{equation}
\end{proposition}

As an immediate consequence of the definition of $d$ and of the local
equivalence of $d_{c}$ with $d$ we obtain the following:

\begin{proposition}
Suppose that $(H_{i}),i=1,2,3$ hold and $X_{0}=x_{0}(\phi ).$ Let $\rho \in
(0,1).$ There exists a universal constant $C$ (depending on $d$ and $\rho $
only) such that for every $R\in (0,1)$ one has%
\begin{equation*}
P(\sup_{t\leq T}d(x_{t}(\phi ),X_{t})\leq R)\geq \exp (-C\mu ^{9}(\frac{T}{h}%
+\int_{0}^{T}\frac{n_{t}^{6(1+d\rho )}}{\lambda _{t}^{1+2d\rho }}(\frac{1}{R}%
+\left\vert \phi _{t}\right\vert ^{2})dt)).
\end{equation*}%
Moreover there exists a constant $C$ (depending on $d$ and $\rho $ but also
on $x_{t}(\phi )$ and on the coefficients $\sigma _{i}(x_{t}(\phi
)),b(x_{t}(\phi ))$ and on their derivatives up to order three) such that%
\begin{equation}
P(\sup_{t\leq T}d_{c}(x_{t}(\phi ),X_{t})\leq R)\geq \exp (-C\mu ^{9}(\frac{T%
}{h}+\int_{0}^{T}\frac{n_{t}^{6(1+d\rho )}}{\lambda _{t}^{1+2d\alpha }}(%
\frac{1}{R}+\left\vert \phi _{t}\right\vert ^{2})dt)).  \label{Not13}
\end{equation}
\end{proposition}

We finish this section with two simple examples.

\smallskip

\textbf{Example 1}. We consider the two dimensional diffusion process
\begin{equation*}
X_{t}^{1}=x_{1}+W_{t}^{1},\quad
X_{t}^{2}=x_{2}+\int_{0}^{t}X_{s}^{1}dW_{s}^{2}.
\end{equation*}%
Straightforward computations give
\begin{equation*}
\left\vert \xi \right\vert _{A_{\delta }(x)}^{2}=\left\vert T_{x,\delta }\xi
\right\vert ^{2}\quad with\quad T_{x,\delta }\xi =(\frac{1}{\sqrt{\delta }}%
\xi _{1},\frac{1}{\sqrt{\delta (\delta +x_{1}^{2})}}\xi _{2}).
\end{equation*}%
In particular, if $x_{1}=0$ then $T_{0,\delta }\xi =(\frac{1}{\sqrt{\delta }}%
\xi _{1},\frac{1}{\delta }\xi _{2})$ and consequently $\{\xi :\left\vert \xi
\right\vert _{A_{\delta }(x)}\leq 1\}$ is an ellipsoid. But if $x_{1}\neq 0$
and $\delta $ is small, then the distance given by $\left\vert \xi
\right\vert _{A_{\delta }(x)}$ is equivalent with the Euclidian one.

If we take a path $x_{t}$ which keeps far from zero then we have ellipticity
along the path and so we may use estimates for elliptic processes (see \cite{[BFM]}
and \cite{[BM]}). But if $x_{1}(t)=0$ for some $t\in \lbrack 0,T]$ then we may no
more use them. Let us compare the norm here and the norm in the elliptic
case: if $x_{1}>0$ the diffusion matrix is not degenerated so we may
consider the norm $\left\vert \xi \right\vert _{B_{\delta }(x)}$ with $%
B_{\delta }(x)=\delta \sigma \sigma ^{\ast }(x).$ We have%
\begin{equation*}
\left\vert \xi \right\vert _{B_{\delta }(x)}^{2}=\frac{1}{\delta }\xi
_{1}^{2}+\frac{1}{\delta x_{1}^{2}}\xi _{2}^{2}\geq \frac{1}{\delta }\xi
_{1}^{2}+\frac{1}{\delta (\delta +x_{1})}\xi _{2}^{2}=\left\vert \xi
\right\vert _{A_{\delta }(x)}^{2}.
\end{equation*}%
So the estimates obtained using the Lie brackets are sharper even if
ellipticity holds.

Let us now take $x_{1}=x_{2}=0,$ $x_{t}(\phi )=(0,0).$ We have $n_{s}=1$ and
$\lambda _{s}=1$ and $X_{t}-x_{t}=(W_{t}^{1},%
\int_{0}^{t}W_{s}^{1}dW_{s}^{2}).$ And we obtain%
\begin{equation*}
P(\sup_{t\leq T}\left( \frac{1}{\delta }\left\vert W_{t}^{1}\right\vert ^{2}+%
\frac{1}{\delta ^{2}}\left\vert \int_{0}^{t}W_{s}^{1}dW_{s}^{2}\right\vert
^{2}\right) \leq 1)=P(\sup_{t\leq T}(\left\vert X_{t}-x_{t}\right\vert
_{A_{\delta }(0)}^{2}\leq 1)\geq e^{-C/\delta }.
\end{equation*}

\textbf{Example 2. The principal invariant diffusion on the Heisenberg group.} %L\'{e}vy area}. 
We consider the diffusion process
\begin{equation*}
X_{t}^{1}=x_{1}+W_{t}^{1},\quad X_{t}^{2}=x_{2}+W_{t}^{2},\quad
X_{t}^{3}=x_{3}+\frac{1}{2}\int_{0}^{t}X_{s}^{1}dW_{s}^{2}-\frac{1}{2}%
\int_{0}^{t}X_{s}^{2}dW_{s}^{1}.
\end{equation*}

Direct computations give%
\begin{equation*}
\left\vert \xi \right\vert _{A_{\delta }(x)}^{2}=\left\vert A_{\delta
}^{-1}(x)\xi \right\vert ^{2}=\frac{1}{\delta }\left( \xi _{1}-\xi
_{3}\times \frac{x_{2}}{2\sqrt{\delta }}\right) ^{2}+\frac{1}{\delta }\left(
\xi _{2}-\xi _{3}\times \frac{x_{1}}{2\sqrt{\delta }}\right) ^{2}+\frac{\xi
_{3}^{2}}{\delta ^{2}}.
\end{equation*}%
In particular for $x=0$ we obtain%
\begin{align*}
&P\Big(\sup_{t\leq T/\delta }\Big( \left\vert W_{t}^{1}\right\vert
^{2}+\left\vert W_{t}^{2}\right\vert ^{2}+A_{t}^{2}(W)\Big) \leq
1\Big)\\
&\qquad=P\Big(\sup_{t\leq T}\Big( \frac{1}{\delta }\left\vert W_{t}^{1}\right\vert
^{2}+\frac{1}{\delta }\left\vert W_{t}^{2}\right\vert ^{2}+\frac{1}{\delta
^{2}}A_{t}^{2}(W)\Big) \leq 1\Big)\geq e^{-\frac{CT}{\delta }}
\end{align*}
where $A_{t}(W)=\int_{0}^{t}W_{s}^{1}dW_{s}^{2}-%
\int_{0}^{t}W_{s}^{2}dW_{s}^{1}.$

\section{Multiple stochastic integrals}

\subsection{Decomposition}

We consider the stochastic process%
\begin{equation}
Z(t)=\sum_{i=1}^{d}a_{i}W_{t}^{i}+\sum_{i,j=1}^{d}a_{i,j}%
\int_{0}^{t}W_{s}^{i}\circ dW_{s}^{j}  \label{Decomp0}
\end{equation}%
with $a_{i},a_{i,j}\in R^{n}.$ Our aim is to give a decomposition for this
process. In order to do it we have to introduce some notation. We fix $%
\delta >0$ and we denote $s_{k}(\delta )=\frac{k}{d}\delta $ and
\begin{equation*}
\Delta _{k}^{i}(\delta ,W)=W_{s_{k}(\delta )}^{i}-W_{s_{k-1}(\delta
)}^{i},\quad \Delta _{k}^{i,j}(\delta ,W)=\int_{s_{k-1}(\delta
)}^{s_{k}(\delta )}(W_{s}^{i}-W_{s_{k-1}}^{i})\circ dW_{s}^{j}.
\end{equation*}%
Notice that $\Delta _{k}^{i,j}(\delta ,W)$ is the Stratonovich integral, but
for $i\neq j$\ it coincides with the Ito integral. When now confusion is
possible we use the short notation $s_{k}=s_{k}(\delta ),\Delta
_{k}^{i}=\Delta _{k}^{i}(\delta ,W),\Delta _{k}^{i,j}=\Delta
_{k}^{i,j}(\delta ,W).$ Moreover for $p=1,...,d$ we define
\begin{equation}\label{Decomp2}
\begin{array}{rcl}
\mu _{p}(\delta ,W)
&=&\displaystyle\sum_{i\neq p}\Delta _{i}^{p}  \\
\psi _{p}(\delta ,W)
&=&\displaystyle\sum_{i\neq j,i\neq p,j\neq p}a_{i,j}\Delta
_{p}^{i,j}+\sum_{l=p+1}^{d}\sum_{i\neq p}\sum_{j\neq l}^{d}a_{i,j}\Delta
_{l}^{j}\Delta _{p}^{i}+\frac{1}{2}\sum_{i\neq p}^{d}a_{i,i}\left\vert
\Delta _{p}^{i}\right\vert ^{2} \\
\varepsilon _{p}(\delta ,W)
&=&\displaystyle\sum_{l>p}^{d}\sum_{j\neq l}^{d}a_{p,j}\Delta
_{l}^{j}+\sum_{p>l}^{d}\sum_{j\neq l}^{d}a_{j,p}\Delta _{l}^{j}+\sum_{j\neq
p}a_{p,j}\Delta _{p}^{j} \\
\eta _{p}(\delta ,W)
&=&\displaystyle\frac{1}{2}a_{p,p}\left\vert \Delta
_{p}^{p}\right\vert ^{2}+\sum_{l>p}^{d}a_{p,l}\Delta _{l}^{l}\Delta
_{p}^{p}+\Delta _{p}^{p}\varepsilon _{p}.
\end{array}
\end{equation}
We denote $\eta (\delta ,W)=\sum_{p=1}^{d}\eta _{p}(\delta ,W)$ and $\psi
(\delta ,W)=\sum_{p=1}^{d}\psi _{p}(\delta ,W)$ and
\begin{equation}
\lbrack a]_{i,p}=a_{i,p}-a_{p,i}.  \label{Decomp3}
\end{equation}%
Our aim is to prove the following decomposition.

\begin{proposition}
\begin{equation}
Z(\delta )=\sum_{p=1}^{d}a_{p}(\Delta _{p}^{p}(\delta ,W)+\mu _{p}(\delta
,W))+\sum_{p=1}^{d}\sum_{i\neq p}[a]_{i,p}\Delta _{p}^{i,p}(\delta ,W)+\eta
(\delta ,W)+\psi (\delta ,W)  \label{Decomp1}
\end{equation}
\end{proposition}

\begin{remark}
The reason of being of this decomposition is the following. We split the
time interval $(0,\delta )$ in $d$ sub intervals of length $\delta /d.$ And
we also split the Brownian motion in corresponding pieces: $%
(W_{s}^{i}-W_{s_{p-1}}^{i})_{s_{p-1}\leq s\leq s_{p}},i=1,...,d.$ Let us fix
$i.$ For $s\in (s_{i-1},s_{i})$ we have the processes $%
(W_{s}^{j}-W_{s_{i-1}}^{j})_{s_{i-1}\leq s\leq s_{i}},j=1,...,d.$ Our idea
is to settle a calculus which is based on $W^{i}$ and to take conditional
expectation with respect to $W^{j},j\neq i.$ So $%
(W_{s}^{j}-W_{s_{i-1}}^{j})_{s_{i-1}\leq s\leq s_{i}},j\neq i$ will appear
as parameters (or controls) which we may choose in an appropriate way. And
the random variables on which the calculus is based are $\Delta
_{i}^{i}=W_{s_{i}}^{i}-W_{s_{i-1}}^{i}$ and $\Delta
_{i}^{j,i}=\int_{s_{i-1}}^{s_{i}}(W_{s}^{j}-W_{s_{i-1}}^{j})dW_{s}^{i},j\neq
i.$ These are the random variables that we have emphasized in the
decomposition of $Z(\delta ).$ Notice that, conditionally to the controls $%
(W_{s}^{j}-W_{s_{i-1}}^{j})_{s_{i-1}\leq s\leq s_{i}},j\neq i,$ this is a
centered Gaussian vector and, under appropriate hypothesis on the controls
this Gaussian vector is non degenerated (we treat in the Appendix 3 the
problem of the choice of the controls). But there is another term which
appear and which is difficult to handle by a choice of the controls $W^{j}:$
this is $\Delta
_{i}^{i,j}=\int_{s_{i-1}}^{s_{i}}(W_{s}^{i}-W_{s_{i-1}}^{i})dW_{s}^{j}.$ So
we use the identity $\Delta _{i}^{i,j}=\Delta _{i}^{j}\Delta _{i}^{i}-\Delta
_{i}^{j,i}$ in order to eliminate this term - and this is the reason for
which $(a_{i,j}-a_{j,i})=[a]_{i,j}$ appears.
\end{remark}

\textbf{Proof}. We decompose
\begin{equation*}
Z(\delta )=\sum_{l=1}^{d}Z(s_{l})-Z(s_{l-1})=\sum_{l=1}^{d}\left(
\sum_{i=1}^{d}a_{i}\Delta
_{l}^{i}+\sum_{i,j=1}^{d}a_{i,j}\int_{s_{l-1}}^{s_{l}}W_{s}^{i}\circ
dW_{s}^{j}\right)
\end{equation*}%
and we write
\begin{equation*}
\int_{s_{l-1}}^{s_{l}}W_{s}^{i}\circ dW_{s}^{j}=W_{s_{l-1}}^{i}\Delta
_{l}^{j}+\Delta _{l}^{i,j}=(\sum_{p=1}^{l-1}\Delta _{p}^{i})\Delta
_{l}^{j}+\Delta _{l}^{i,j}.
\end{equation*}%
Then%
\begin{equation*}
Z(\delta )=\sum_{l=1}^{d}\sum_{i=1}^{d}a_{i}\Delta
_{l}^{i}+\sum_{l=1}^{d}\sum_{i,j=1}^{d}a_{i,j}(\sum_{p=1}^{l-1}\Delta
_{p}^{i})\Delta _{l}^{j}+\sum_{l=1}^{d}\sum_{i,j=1}^{d}a_{i,j}\Delta
_{l}^{i,j}=:S_{1}+S_{2}+S_{3}.
\end{equation*}%
Notice first that%
\begin{equation*}
S_{1}=\sum_{l=1}^{d}a_{l}\Delta _{l}^{l}+\sum_{l=1}^{d}\sum_{i\neq
l}a_{i}\Delta _{l}^{i}.
\end{equation*}%
\bigskip We treat now $S_{3}.$ We will use the identities%
\begin{equation*}
\left\vert \Delta _{l}^{i}\right\vert ^{2}=2\Delta _{l}^{i,i}\quad and\quad
\Delta _{l}^{i}\Delta _{l}^{j}=\Delta _{l}^{i,j}+\Delta _{l}^{j,i}.
\end{equation*}%
Then%
\begin{eqnarray*}
.S_{3} &=&\sum_{l=1}^{d}\sum_{i=1}^{d}a_{i,i}\Delta
_{l}^{i,i}+\sum_{l=1}^{d}\sum_{i\neq j}a_{i,j}\Delta _{l}^{i,j} \\
&=&\sum_{l=1}^{d}\sum_{i=1}^{d}a_{i,i}\Delta
_{l}^{i,i}+\sum_{l=1}^{d}\sum_{i\neq l}a_{i,l}\Delta
_{l}^{i,l}+\sum_{l=1}^{d}\sum_{j\neq l}a_{l,j}\Delta
_{l}^{l,j}+\sum_{l=1}^{d}\sum_{i\neq j,i\neq lj\neq l}a_{i,j}\Delta
_{l}^{i,j} \\
&=&\frac{1}{2}\sum_{l=1}^{d}\sum_{i=1}^{d}a_{i,i}\left\vert \Delta
_{l}^{i}\right\vert ^{2}+\sum_{l=1}^{d}\sum_{i\neq l}a_{i,l}\Delta _{l}^{i,l}
\\
&&+\sum_{l=1}^{d}\sum_{j\neq l}a_{l,j}\left( \Delta _{l}^{j}\Delta
_{l}^{l}-\Delta _{l}^{j,l}\right) +\sum_{l=1}^{d}\sum_{i\neq j,i\neq l,j\neq
l}a_{i,j}\Delta _{l}^{i,j} \\
&=&\frac{1}{2}\sum_{i=1}^{d}a_{i,i}\left\vert \Delta _{i}^{i}\right\vert
^{2}+\frac{1}{2}\sum_{l=1}^{d}\sum_{i\neq l}^{d}a_{i,i}\left\vert \Delta
_{l}^{i}\right\vert ^{2}+\sum_{l=1}^{d}\sum_{i\neq l}(a_{i,l}-a_{l,i})\Delta
_{l}^{i,l} \\
&&+\sum_{l=1}^{d}\left( \sum_{j\neq l}a_{l,j}\Delta _{l}^{j}\right) \Delta
_{l}^{l}+\sum_{l=1}^{d}\sum_{i\neq j,i\neq l,\neq j\neq }a_{i,j}\Delta
_{l}^{i,j}.
\end{eqnarray*}%
We treat now $S_{2}.$ We want to emphasis terms which contain $\Delta
_{i}^{i}.$ We have
\begin{equation*}
S_{2}=\sum_{l>p}^{d}\sum_{i,j=1}^{d}a_{i,j}\Delta _{p}^{i}\Delta
_{l}^{j}=S_{2}^{\prime }+S_{2}^{\prime \prime }+S_{2}^{\prime \prime \prime
}+S_{2}^{iv}
\end{equation*}%
with $\sum_{l>p}^{d}=\sum_{p=1}^{d}\sum_{l=p+1}^{d}$ and%
\begin{eqnarray*}
S_{2}^{\prime } &=&\sum_{l>p}^{d}a_{p,l}\Delta _{p}^{p}\Delta _{l}^{l},\quad
S_{2}^{\prime \prime }=\sum_{l>p}^{d}\sum_{j\neq l}^{d}a_{p,j}\Delta
_{p}^{p}\Delta _{l}^{j} \\
S_{2}^{\prime \prime \prime } &=&\sum_{l>p}^{d}\sum_{i\neq
p}^{d}a_{i,l}\Delta _{p}^{i}\Delta _{l}^{l},\quad
S_{2}^{iv}=\sum_{l>p}^{d}\sum_{i\neq p,j\neq l}^{d}a_{i,j}\Delta
_{p}^{i}\Delta _{l}^{j}.
\end{eqnarray*}%
We have
\begin{equation*}
S_{2}^{\prime \prime }=\sum_{p=1}^{d}\Delta _{p}^{p}\left(
\sum_{l=p+1}^{d}\sum_{j\neq l}^{d}a_{p,j}\Delta _{l}^{j}\right)
\end{equation*}%
and
\begin{equation*}
S_{2}^{\prime \prime \prime }=\sum_{l=1}^{d}\Delta _{l}^{l}\left(
\sum_{p=1}^{l-1}\sum_{i\neq p}^{d}a_{i,l}\Delta _{p}^{i}\right)
=\sum_{p=1}^{d}\Delta _{p}^{p}\left( \sum_{l=1}^{p-1}\sum_{j\neq
l}^{d}a_{j,p}\Delta _{l}^{j}\right)
\end{equation*}%
so that%
\begin{equation*}
S_{2}^{\prime \prime }+S_{2}^{\prime \prime \prime }=\sum_{p=1}^{d}\Delta
_{p}^{p}\left( \sum_{l=p+1}^{d}\sum_{j\neq l}^{d}a_{p,j}\Delta
_{l}^{j}+\sum_{l=1}^{p-1}\sum_{j\neq l}^{d}a_{j,p}\Delta _{l}^{j}\right) .
\end{equation*}

Finally%
\begin{eqnarray*}
Z(\delta ) &=&\sum_{l=1}^{d}a_{l}\Delta _{l}^{l}+\sum_{l=1}^{d}\sum_{i\neq
l}a_{i}\Delta _{l}^{i} \\
&&+\sum_{l>p}^{d}a_{p,l}\Delta _{p}^{p}\Delta _{l}^{l}+\sum_{p=1}^{d}\Delta
_{p}^{p}\left( \sum_{l>p}^{d}\sum_{j\neq l}^{d}a_{p,j}\Delta
_{l}^{j}+\sum_{p>l}^{d}\sum_{j\neq l}^{d}a_{j,p}\Delta _{l}^{j}\right)\\
&&+\sum_{l>p}^{d}\sum_{i\neq p,j\neq l}^{d}a_{i,j}\Delta _{p}^{i}\Delta
_{l}^{j}
+\frac{1}{2}\sum_{i=1}^{d}a_{i,i}\left\vert \Delta _{i}^{i}\right\vert
^{2}+\frac{1}{2}\sum_{l=1}^{d}\sum_{i\neq l}^{d}a_{i,i}\left\vert \Delta
_{l}^{i}\right\vert ^{2}\\
&&+\sum_{l=1}^{d}\sum_{i\neq l}(a_{i,l}-a_{l,i})\Delta
_{l}^{i,l}
+\sum_{l=1}^{d}\left( \sum_{j\neq l}a_{l,j}\Delta _{l}^{j}\right) \Delta
_{l}^{l}+\sum_{l=1}^{d}\sum_{i\neq j,i\neq l,j\neq l}a_{i,j}\Delta
_{l}^{i,j}.
\end{eqnarray*}%
We want to compute the coefficient of $\Delta _{p}^{p}:$ this term appears
in
\begin{eqnarray*}
&&\sum_{p=1}^{d}\Delta _{p}^{p}(a_{p}+\varepsilon _{p})\quad with \\
\varepsilon _{p} &=&\sum_{l>p}^{d}\sum_{j\neq l}^{d}a_{p,j}\Delta
_{l}^{j}+\sum_{p>l}^{d}\sum_{j\neq l}^{d}a_{j,p}\Delta _{l}^{j}+\sum_{j\neq
p}a_{p,j}\Delta _{p}^{j}.
\end{eqnarray*}%
We consider now $\Delta _{p}^{i,p}.$ It appears in
\begin{equation*}
\sum_{p=1}^{d}\sum_{i\neq p}(a_{i,p}-a_{p,i})\Delta _{p}^{i,p}
\end{equation*}%
The other terms are%
\begin{eqnarray*}
&&\sum_{l=1}^{d}\sum_{i\neq l}a_{i}\Delta _{l}^{i}+\sum_{l>p}^{d}\sum_{i\neq
p,j\neq l}^{d}a_{i,j}\Delta _{p}^{i}\Delta _{l}^{j}+\frac{1}{2}%
\sum_{i=1}^{d}a_{i,i}\left\vert \Delta _{i}^{i}\right\vert ^{2}+\frac{1}{2}%
\sum_{l=1}^{d}\sum_{i\neq l}^{d}a_{i,i}\left\vert \Delta _{l}^{i}\right\vert
^{2} \\
&&+\sum_{l=1}^{d}\sum_{i\neq j,i\neq l,j\neq l}a_{i,j}\Delta
_{l}^{i,j}+\sum_{l=p+1}^{d}a_{p,l}\Delta _{p}^{p}\Delta _{l}^{l}.
\end{eqnarray*}%
We put everything together and (\ref{Decomp1}) is proved. $\square $

\subsection{Main estimates}

Throughout this section we will assume that%
\begin{equation}
Span\{a_{i},[a]_{j,p},i,j,p=1,...,d,j\neq p\}=R^{n}.  \label{Estimate0}
\end{equation}%
Let us introduce some notation. We consider the matrix $%
A=(a_{i},[a]_{j,p},i,j,p=1,...,d,j\neq p)$ to be the matrix with columns $%
a_{i}$ and $[a]_{j,p}.$ For $R\in (0,1]$ we define the matrix $A_{R}=(\sqrt{R%
}a_{i},R[a]_{j,p},i,j,p=1,...,d,j\neq p)$ and we denote $\lambda _{\ast
}(A_{R}),\lambda ^{\ast }(A_{R})$ the lower and the larger eigenvalue of $%
A_{R}A_{R}^{\ast }.$ We just write $\lambda _{\ast }(A),\lambda ^{\ast }(A)$
if $R=1.$ We associate the norms $\left\vert y\right\vert
_{A_{R}}^{2}=\left\langle (A_{R}A_{R})^{-1}y,y\right\rangle .$

In Proposition \ref{NORM1} from the Appendix 4 we prove the following
basic properties. For every $0<R\leq R^{\prime }\leq 1$%
\begin{equation}
\sqrt{\frac{R}{R^{\prime }}}\left\vert y\right\vert _{A_{R}}\geq \left\vert
y\right\vert _{A_{R^{\prime }}}\geq \frac{R}{R^{\prime }}\left\vert
y\right\vert _{A_{R}}  \label{Estimate1}
\end{equation}%
and
\begin{equation}
\frac{1}{\sqrt{R}\sqrt{\lambda ^{\ast }(A)}}\left\vert y\right\vert \leq
\left\vert y\right\vert _{A_{R}}\leq \frac{1}{R\sqrt{\lambda _{\ast }(A)}}%
\left\vert y\right\vert .  \label{Estimate2}
\end{equation}%
Finally
\begin{equation}
\left\vert A_{R}y\right\vert _{A_{R}}\leq \left\vert y\right\vert .
\label{Estimate2'}
\end{equation}

\begin{lemma}\label{ESTIMATE1}
Suppose that (\ref{Estimate0}) holds. There
exists an universal constant $C_{0}$ such that for every $R\geq \delta >0$
and $r>0$%
\begin{equation}
P(\sup_{t\leq \delta }\left\vert Z_{t}\right\vert _{A_{R}}\geq r)\leq \exp \Big(-%
\frac{rR}{C_{0}\delta }\Big( r\wedge \frac{\sqrt{\lambda _{\ast }(A)}}{%
\overline{a}}\Big) \Big)  \label{Estimate8}
\end{equation}%
with%
\begin{equation}
\overline{a}=1\vee \max_{i,j}\left\vert a_{i,j}\right\vert .
\label{Estimate9}
\end{equation}%
\end{lemma}

\begin{remark}
One might think to use directly Bernstein's inequality in order to estimate $%
P(\sup_{t\leq \delta }\left\vert Z_{t}\right\vert _{A_{R}}\geq r)$ but then
one would  not obtain the right inequality. Indeed one writes $\left\vert
Z_{t}\right\vert _{A_{R}}\leq (R\sqrt{\lambda _{\ast }(A)})^{-1}\left\vert
Z_{t}\right\vert $ and then the above probability is bounded by%
\begin{equation*}
P(\sup_{t\leq \delta }\left\vert Z_{t}\right\vert \geq rR\sqrt{\lambda
_{\ast }(A)})\leq \exp (-\frac{r^{2}R^{2}\lambda _{\ast }(A)}{\delta }).
\end{equation*}%
So one obtains $\frac{R^{2}}{\delta }$ instead of $\frac{R}{\delta }$ and
this is not in the right scale. The reason is that in the above argument we
just use the lower eigenvalue $\lambda _{\ast }(A)$ in order to upper bound $%
\left\vert Z_{t}\right\vert _{A_{R}}$ since in the proof of our lemma we use
the more subtle inequality $\left\vert A_{R}y\right\vert _{A_{R}}\leq
\left\vert y\right\vert .$
\end{remark}

\textbf{Proof.} Let $t\leq \delta$. We decompose $Z(t)$ instead of $Z(\delta)$ and similarly to (\ref{Decomp1}) we obtain
\begin{equation*}
Z(t)=\sum_{p=1}^{d}a_{p}(\Delta _{p}^{p}(t,W)+\mu
_{p}(t,W))+\sum_{p=1}^{d}\sum_{i\neq p}[a]_{i,p}\Delta _{p}^{i,p}(t,W)+\eta
(t,W)+\psi (t,W),
\end{equation*}%
in which $\eta(t,W)$ and $\psi (t,W)$ are defined as in (\ref{Decomp2}) with $\Delta^i_p$ and $\Delta^{ij}_p$ replaced by $\Delta^{i}_{p}(t,W)$ and $\Delta^{ij}_{p}(t,W)$ respectively, and these last quantities are defined as follows: for $t\in [0,T]$,
$\Delta^i_p(t,W)=W^i_{s_p\wedge t}-W^i_{s_{p-1}\wedge t}$ and $\Delta^{ij}_p(t,W)=\int_{s_{p-1}\wedge t}^{s_{p}\wedge t}(W^i_s-W^i_{s_{p-1}\wedge t})dW^j_s$.

We denote by $u(t)\in R^{m}$ the vector with component $u_{p}(t)=t^{-1/2}(%
\Delta _{p}^{p}(t,W)+\mu _{p}(t,W))=t^{-1/2}W_{t}^{p},p=1,...,d$ and $%
u_{i,j}(t)=0,i\neq j$ and we also denote%
\begin{equation*}
U(t)=\sum_{p=1}^{d}\sum_{i\neq p}[a]_{i,p}\Delta _{p}^{i,p}(t,W)+\eta
(t,W)+\psi (t,W).
\end{equation*}%
Then we have%
\begin{equation*}
Z(t)=\sum_{p=1}^{d}t^{1/2}a_{p}u_{p}(t)+\sum_{p=1}^{d}\sum_{i\neq
p}t[a]_{i,p}\times 0+U(t)=A_{t}u(t)+U(t).
\end{equation*}%
Using the norm inequalities given above%
\begin{equation*}
\left\vert U(t)\right\vert _{A_{R}}\leq \frac{1}{R\sqrt{\lambda _{\ast }(A)}}%
\left\vert U(t)\right\vert \leq \frac{C\overline{a}}{R\sqrt{\lambda _{\ast
}(A)}}\sum_{i,j=1}^{d}(\left\vert \Delta _{j}^{i}(t,W)\right\vert
^{2}+\sum_{p=1}^{d}\left\vert \Delta _{p}^{i,j}(t,W)\right\vert )
\end{equation*}%
so that%
\begin{align*}
P\Big(\sup_{t\leq \delta }\left\vert U(t)\right\vert _{A_{R}} \geq \frac{r}{2}%
\Big)\leq &\sum_{i,j=1}^{d}P\Big(\sup_{t\leq \delta }\left\vert \Delta
_{i}^{j}(t,W)\right\vert ^{2}\geq \frac{rR\sqrt{\lambda _{\ast }(A)}}{C%
\overline{a}}\Big) \\
&+\sum_{i,j,p=1}^{d}P\Big(\sup_{t\leq \delta }\left\vert \Delta
_{p}^{i,j}(t,W)\right\vert \geq \frac{rR\sqrt{\lambda _{\ast }(A)}}{C%
\overline{a}}\Big).
\end{align*}%
It is easy to check that
\begin{equation*}
P\Big(\sup_{t\leq \delta }\left\vert \Delta _{p}^{p}(t,W)\right\vert ^{2}\geq
\frac{rR\sqrt{\lambda _{\ast }(A)}}{C\overline{a}}\Big)\leq C^{\prime }\exp \Big(-%
\frac{rR\sqrt{\lambda _{\ast }(A)}}{C^{\prime }\overline{a}\delta }\Big).
\end{equation*}%
Moreover,
\begin{equation*}
\sup_{t\leq \delta }\left\vert \Delta _{p}^{i,j}(t,W)\right\vert \leq
2\sup_{t\leq \delta }\left\vert \int_{0}^{t}W_{s}^{i}dW_{s}^{j}\right\vert
+2\sup_{t\leq \delta }(\left\vert W_{t}^{i}\right\vert ^{2}+\left\vert
W_{t}^{j}\right\vert ^{2}).
\end{equation*}%
Using (\ref{Exp1}) from the Appendix 1 we obtain%
\begin{equation*}
P\Big(\sup_{t\leq \delta }\Big\vert \int_{0}^{t}W_{s}^{i}dW_{s}^{j}\Big\vert
\geq \frac{rR\sqrt{\lambda _{\ast }(A)}}{C\overline{a}}\Big)\leq C\exp \Big(-\frac{rR%
\sqrt{\lambda _{\ast }(A)}}{C\overline{a}\delta }\Big).
\end{equation*}%
So we have proved that
\begin{equation*}
P\Big(\sup_{t\leq \delta }\left\vert U(t)\right\vert _{A_{R}}\geq \frac{r}{2}%
\Big)\leq C\exp \Big(-\frac{rR\sqrt{\lambda _{\ast }(A)}}{C\overline{a}\delta }\Big).
\end{equation*}%
Using (\ref{Estimate1}) (recall that $t\leq \delta \leq R)$ and (\ref%
{Estimate2'})
\begin{equation*}
\left\vert A_{t}u(t)\right\vert _{A_{R}}\leq \sqrt{\frac{t}{R}}\left\vert
A_{t}u(t)\right\vert _{A_{t}}\leq \sqrt{\frac{t}{R}}\left\vert
u(t)\right\vert \leq \frac{C}{\sqrt{R}}\sup_{t\leq \delta }\left\vert
W_{t}\right\vert .
\end{equation*}%
It follows that
\begin{equation*}
P\Big(\sup_{t\leq \delta }\left\vert A_{t}u(t)\right\vert _{A_{R}}\geq \frac{r}{2%
}\Big)\leq P\Big(\sup_{t\leq \delta }\left\vert W_{t}\right\vert \geq \frac{r\sqrt{R}%
}{C}\Big)\leq C\exp \Big(-\frac{r^{2}R}{C\delta }\Big).
\end{equation*}%
$\square $

\medskip

We give the main result in this section.

\begin{proposition}
\label{ESTIMATE2}Suppose that $\lambda _{\ast }(A)>0.$ Let $\rho \in (0,1)$
be fixed. There exists an universal constant $C_{\ast }$ (depending on $d$
and on $\rho $ only) such that for every
\begin{equation}
r\leq \frac{\lambda _{\ast }^{1/2}(A)}{C_{\ast }\overline{a}}
\label{Estimate6'}
\end{equation}
one has%
\begin{equation}
P(\left\vert Z_{\delta }\right\vert _{A_{\delta }}\leq r)\geq \frac{r^{m}}{%
C_{\ast }}\times \frac{\lambda _{\ast }^{2d^{3}}(A)}{\overline{a}^{d^{3}}}%
\times \exp (-\frac{C_{\ast }\lambda _{\ast }^{d^{2}\rho }(A)}{\overline{a}%
^{2}}).  \label{Estimate6bis}
\end{equation}
\end{proposition}

\textbf{Proof. Step 1}. \textbf{Scaling}. Let $B_{t}=\delta
^{-1/2}W_{t\delta }.$ Then $B$ is a standard Brownian motion and we denote%
\begin{equation*}
\Delta _{i}^{j}(B)=B_{i}^{j}-B_{i-1}^{j},\quad \Delta
_{p}^{i,j}(B)=\int_{p-1}^{p}(B_{s}^{j}-B_{p}^{j})dB_{s}^{i},i\neq j.
\end{equation*}%
We also denote by $\Delta (B)$ the vector $(\Delta _{i}^{j}(B),\Delta
_{p}^{i,j}(B),i,j,p=1,...,d)$ and we define $\Theta (B)=(\Theta
_{1}(B),...,\Theta _{d}(B))$ with $\Theta _{p}(B)=(\Delta _{p}^{p}(B),\Delta
_{p}^{j,p}(B),j\neq p).$ We consider the $\sigma $ field
\begin{equation*}
\mathcal{G}:=\sigma (W_{s}^{j}-W_{s_{p-1}(\delta )}^{j},s_{p-1}(\delta )\leq
s\leq s_{p}(\delta ),p=1,...d,j\neq p).
\end{equation*}%
Conditionally to $\mathcal{G}$ the random variable $\Theta _{p}(B)$ is
Gaussian with covariance matrix $Q_{p}(B)$ given by%
\begin{align*}
Q_{p}^{p,j}(B)& =\int_{p-1}^{p}(B_{s}^{j}-B_{i-1}^{j})ds,\quad j\neq
p, \\
Q_{p}^{i,j}(B)&
=\int_{p-1}^{p}(B_{s}^{j}-B_{i-1}^{j})(B_{s}^{i}-B_{i-1}^{i})ds,\quad j\neq
p,i\neq p, \\
Q_{p}^{p,p}(B)& =1.
\end{align*}%
Since the random variables $\Theta _{1}(B),...,\Theta _{d}(B)$ are
independent $\Theta (B)$ is a Gaussian random variable. We denote by $Q(B)$
the covariance matrix of $\Theta (B)$ and by $\lambda _{\ast }(B),\lambda
^{\ast }(B)$ the smaller and the larger eigenvalues of $Q(B).$ Since this
matrix is built with the blocks $Q_{p}(B),p=1,...,d$ we have%
\begin{equation*}
\lambda _{\ast }(B)=\prod_{p=1}^{d}\lambda _{\ast ,p}(B)\quad and\quad
\lambda _{\ast }(B)=\prod_{p=1}^{d}\lambda _{p}^{\ast }(B)
\end{equation*}%
where $\lambda _{\ast ,p}(B),\lambda _{p}^{\ast }(B)$ are the smaller and
the larger eigenvalues of $Q_{p}(B).$

We come now back to our problem. Let $\eta (\Delta (B)),\psi (\Delta
(B)),\varepsilon (\Delta (B)),\mu (\Delta (B))$ be the quantities defined in
(\ref{Decomp2}) with $\Delta =\Delta (\delta ,W)$ replaced by $\Delta (B).$
Then $\delta \eta (\Delta (B))=\eta (\delta ,W).$ The same is true for $\psi
$ and $\varepsilon $ and finally $\sqrt{\delta }\mu (\Delta (B))=\mu (\delta
,W).$ So using (\ref{Decomp1})%
\begin{equation*}
Z_{\delta }=\sum_{p=1}^{d}\sqrt{\delta }a_{p}(\Delta _{p}^{p}(B)+\mu
_{p}(\Delta (B)))+\sum_{p=1}^{d}\sum_{i\neq p}\delta \lbrack a]_{ip}\Delta
_{p}^{i,p}(B)+\delta \eta (\Delta )+\delta \psi (\Delta ).
\end{equation*}%
We define now the vector $\mu (\Delta (B))=(\mu _{p}(\Delta (B)),\mu
_{i,j}(\Delta (B)\in R^{m},i\neq j)$ by $\mu _{i,j}(\Delta (B))=0$ and then
we may write the above decomposition in matrix notation%
\begin{eqnarray}
Z_{\delta } &=&A_{\delta }(\Theta (B)+\mu (\Delta (B)))+\delta \eta (\Theta
(B))+\delta \psi (\Delta (B))  \label{Decomp12} \\
&=&y+.A_{\delta }\Theta (B)+\eta _{\delta }(\Theta (B))  \notag
\end{eqnarray}%
with%
\begin{equation*}
y=A_{\delta }\mu (\Delta (B))+\delta \psi (\Delta (B)),\quad \eta _{\delta
}(\theta )=\delta \eta (\theta ).
\end{equation*}

\textbf{Step 2. Localization}. We take%
\begin{equation}
\varepsilon \leq \frac{\lambda _{\ast }(A)}{C_{1}\overline{a}^{2}}
\label{Estimate4''}
\end{equation}%
where $C_{1}$ is an universal constant to be chosen in the sequel. For each
$p=1,...,d$ we define the sets%
\begin{equation*}
\Lambda _{\rho ,\varepsilon ,p}=\Big\{\det Q_{p}(B)\geq \varepsilon ^{\rho
},\sup_{p-1\leq t\leq p}\sum_{j\neq p}\left\vert
B_{t}^{j}-B_{p-1}^{j}\right\vert \leq \varepsilon ^{-\rho },q_{p}(B)\leq
\varepsilon \Big\}
\end{equation*}%
with%
\begin{equation*}
q_{p}(B)=\sum_{j\neq p}\left\vert B_{p}^{j}-B_{p-1}^{j}\right\vert
+\sum_{j\neq p,i\neq p}\left\vert
\int_{p-1}^{p}(B_{s}^{j}-B_{i-1}^{j})dB_{s}^{i}\right\vert .
\end{equation*}%
By (\ref{Sup2}) in Appendix 3 we may find some constants $c$ and $%
\varepsilon _{\ast }$ depending on $d$ and $\rho $ only such that
\begin{equation}
P(\Lambda _{\rho ,\varepsilon ,p})\geq c\varepsilon ^{\frac{1}{2}%
d(d+1)}\quad \mbox{for}\quad \varepsilon \leq \varepsilon _{\ast }
\label{Estimate5'}
\end{equation}%
And using the independence we obtain%
\begin{equation}
P\big(\cap _{p=1}^{d}\Lambda _{\rho ,\varepsilon ,p}\big)\geq c^{d}\times
\varepsilon ^{\frac{1}{2}d^{2}(d+1)}.  \label{Estimate4}
\end{equation}%
On the set $\cap _{p=1}^{d}\Lambda _{\rho ,\varepsilon ,p}$ we have $\det
Q_{p}(B)\geq \varepsilon ^{\rho }$ so that $\det Q(B)\geq \varepsilon
^{d\rho }.$ We also have $\lambda ^{\ast }(B)\leq \varepsilon ^{-\rho }$ and
this gives $\lambda _{\ast }(B)\geq \varepsilon ^{d^{2}\rho }.$ And we also
have $\det Q(B)\leq \varepsilon ^{-d\rho }$ so
\begin{equation}
\cap _{p=1}^{d}\Lambda _{\rho ,\varepsilon ,p}\subset \big\{\det Q(B)\leq
\varepsilon ^{-d\rho },\lambda _{\ast }(B)\geq \varepsilon ^{d^{2}\rho
},\sum_{p=1}^{d}q_{p}(B)\leq d\varepsilon \big\}  \label{Estimate4'}
\end{equation}

\textbf{Step 3. Inverse function theorem}. We will use (\ref{DENSITY4}) with $%
G=Z_{\delta }$ so we have to estimate the parameters associated to $\eta
_{\delta }$ and $A_{\delta }.$ Notice first that $\lambda _{\ast }(A_{\delta
})\geq \delta ^{2}\lambda _{\ast }(A),c_{3,\eta _{\delta }}=0$ and $%
c_{2,\eta _{\delta }}\leq C_{2}\overline{a}\delta .$ So the first inequality
in (\ref{DENSITY6}) reads%
\begin{equation*}
r\leq \frac{\lambda _{\ast }^{1/2}(A)}{C_{2}\overline{a}}\leq \frac{\lambda
_{\ast }^{1/2}(A_{\delta })}{16(c_{2,\eta _{\delta }}+c_{3,\eta _{\delta }})}%
.
\end{equation*}%
And this is verified by our hypothesis. Moreover%
\begin{equation*}
c_{\ast }(\eta _{\delta },r)\leq C_{3}\overline{a}(\left\vert \theta
\right\vert +\sum_{p=1}^{d}\left\vert \varepsilon _{p}(\Delta
(B))\right\vert )\leq C_{4}\overline{a}(r+\sum_{p=1}^{d}q_{p}(B))\leq C_{4}%
\overline{a}(\frac{\lambda _{\ast }^{1/2}(A)}{C_{2}\overline{a}}%
+d\varepsilon ).
\end{equation*}%
If we choose $C_{1}$ in (\ref{Estimate4''}) sufficiently large and $C_{2}$
large also we obtain $c_{\ast }(\eta _{\delta },r)\leq \frac{1}{2}$ which is
the second restriction in (\ref{DENSITY6}). Let $p_{\mathcal{G},Z_{\delta
}}(z)$ be the density of $Z_{\delta }$ conditionally to $\mathcal{G}.$ Then,
using (\ref{DENSITY4}), if $\left\vert z-y\right\vert _{A_{\delta }}\leq
r\leq 1$\ we obtain
\begin{eqnarray*}
p_{\mathcal{G},Z_{\delta }}(z) &\geq &\frac{(4\lambda _{\ast }(B))^{(m-n)/2}%
}{(8\pi )^{m/2}\sqrt{\det Q(B)}\sqrt{\det A_{\delta }A_{\delta }^{\ast }}}%
\exp (-\frac{1}{4\lambda _{\ast }(Q(B))}\left\vert z-y\right\vert
_{A_{\delta }}^{2}) \\
&\geq &\frac{\varepsilon ^{d^{3\rho }}}{(8\pi )^{m/2}\sqrt{\det A_{\delta
}A_{\delta }^{\ast }}}\exp (-\frac{1}{4\varepsilon ^{d^{2}\rho }})
\end{eqnarray*}%
the second inequality being true on $\cap _{p=1}^{d}\Lambda _{\rho
,\varepsilon ,p}.$ On this set we also have
\begin{equation*}
\left\vert \mu (\Delta (B))\right\vert +\left\vert \psi (\Delta
(B))\right\vert \leq C_{5}\overline{a}\sum_{p=1}^{d}q_{p}(B)\leq C_{6}%
\overline{a}\varepsilon
\end{equation*}%
so that%
\begin{eqnarray*}
\left\vert y\right\vert _{A_{\delta }} &\leq &\left\vert A_{\delta }\mu
(\Delta (B))\right\vert _{A_{\delta }}+\delta \left\vert \psi (\Delta
(B))\right\vert _{A_{\delta }}\leq \left\vert \mu (\Delta (B))\right\vert +%
\frac{1}{\sqrt{\lambda _{\ast }(A)}}\left\vert \psi (\Delta (B))\right\vert
\\
&\leq &\frac{C_{7}\overline{a}}{\sqrt{\lambda _{\ast }(A)}}\varepsilon \leq
\frac{r}{2}.
\end{eqnarray*}%
So, if $\left\vert z\right\vert _{A_{\delta }}\leq \frac{r}{2}$ then $%
\left\vert z-y\right\vert _{A_{\delta }}\leq r.$ It follows that
\begin{eqnarray*}
P_{\mathcal{G}}(\left\vert Z_{\delta }\right\vert _{A_{\delta }} &\leq &%
\frac{r}{2})=\int_{\{\left\vert z\right\vert _{A_{\delta }}\leq \frac{r}{2}%
\}}p_{\mathcal{G},Z_{\delta }}(z)dz\geq \frac{\varepsilon ^{d^{3\rho }}}{%
(8\pi )^{m/2}}\exp (-\frac{1}{4\varepsilon ^{d^{2}\rho }})\int_{\{\left\vert
z\right\vert _{A_{\delta }}\leq \frac{r}{2}\}}\frac{1}{\sqrt{\det A_{\delta
}A_{\delta }^{\ast }}}dz \\
&=&\frac{\varepsilon ^{d^{3\rho }}}{(8\pi )^{m/2}}\exp (-\frac{1}{%
4\varepsilon ^{d^{2}\rho }})\times \frac{r^{m}}{2^{m}}
\end{eqnarray*}%
the last equality being obtained by a change of variable. Finally using (\ref%
{Estimate4})%
\begin{equation*}
P(\left\vert Z_{\delta }\right\vert _{A_{\delta }}\leq \frac{r}{2})\geq P(P_{%
\mathcal{G}}(\left\vert Z_{\delta }\right\vert _{A_{\delta }}\leq r),\cap
_{p=1}^{d}\Lambda _{\rho ,\varepsilon ,p})\geq \frac{r^{m}\varepsilon
^{2d^{3}}}{C_{8}}\exp (-\frac{1}{4\varepsilon ^{d^{2}\rho }}).
\end{equation*}%
We replace now $\varepsilon $ by the expression in the RHS of (\ref%
{Estimate4''}) and we obtain (\ref{Estimate6bis}). $\square $

\begin{corollary}
\label{ESTIMATE3}Suppose that $\lambda _{\ast }(A)>0.$ Let $\rho \in (0,1)$
be fixed. There exists some universal constant $C$ (depending on $d$ and on $%
\rho $ only) such that for every $r,R>0$ the following holds. Suppose that
\begin{equation}
\delta \leq \frac{rR}{C\ln \frac{1}{r}}\left( r\wedge \frac{\sqrt{\lambda
_{\ast }(A)}}{\overline{a}}\right) \times \frac{\lambda _{\ast }^{d\rho }(A)%
}{\overline{a}^{2d\lambda }}.  \label{Estimate13}
\end{equation}%
Then%
\begin{equation}
P(\sup_{t\leq \delta }\left\vert Z_{t}\right\vert _{A_{R}}\leq r,\left\vert
Z_{\delta }\right\vert _{A_{\delta }}\leq r)\geq \frac{r^{m}}{2C_{\ast }}%
\exp (-\frac{C_{\ast }\overline{a}^{2d\rho }}{\lambda _{\ast }^{d\lambda }(A)%
})  \label{Estimate12}
\end{equation}%
with $C_{\ast }$ the constant from (\ref{Estimate6bis}).
\end{corollary}

\textbf{Proof}. We use (\ref{Estimate8}) and (\ref{Estimate6bis}) in order
to obtain
\begin{eqnarray*}
P(\sup_{t\leq \delta }\left\vert Z_{t}\right\vert _{A_{R}} &\leq
&r,\left\vert Z_{\delta }\right\vert _{A_{\delta }}\leq r)\geq P(\left\vert
Z_{\delta }\right\vert _{A_{\delta }}\leq r)-P(\sup_{t\leq \delta
}\left\vert Z_{t}\right\vert _{A_{R}}>r) \\
&\geq &\frac{r^{m}}{C_{3}}\exp (-\frac{C_{3}\overline{a}^{2d\rho }}{\lambda
_{\ast }^{d\lambda }(A)})-\exp (-\frac{rR}{C_{0}\delta }\left( r\wedge \frac{%
\sqrt{\lambda _{\ast }(A)}}{\overline{a}}\right) ) \\
&\geq &\frac{r^{m}}{2C_{3}}\exp (-\frac{C_{3}\overline{a}^{2d\rho }}{\lambda
_{\ast }^{d\lambda }(A)})
\end{eqnarray*}%
the last inequality being a consequence of our restriction on $\delta
.\square $

\section{Diffusion processes}\label{diff-proc}

\subsection{Short time behavior}\label{diff-proc-1}

We consider the diffusion process $X_{t}$ solution of (\ref{Not1}) and the
skeleton $x_{t}=x_{t}(\phi )$ solution of (\ref{Not3}) and we give for them
an estimate which is analogous to (\ref{Estimate12}). Using a development in
stochastic Taylor series of order two we write%
\begin{equation*}
X_{t}=X_{0}+Z_{t}+b(0,X_{0})t+R_{t}
\end{equation*}%
where $Z_{t}$ is defined in (\ref{Decomp0}) with $a_{i}=\sigma
_{i}(0,X_{0})$, $a_{i,j}=\partial _{\sigma _{i}}\sigma _{j}(0,X_{0})$ so
that $[a]_{i,j}=[\sigma _{i},\sigma _{j}](0,X_{0})$, and
\begin{align*}
R_{t}=& \sum_{j,i=1}^{d}\int_{0}^{t}\int_{0}^{s}(\partial _{\sigma
_{i}}\sigma _{j}(u,X_{u})-\partial _{\sigma _{i}}\sigma _{j}(0,X_{0}))\circ
dW_{u}^{i}\circ dW_{s}^{j}\\
& +\sum_{i=1}^{d}\int_{0}^{t}\int_{0}^{s}\partial
_{b}\sigma _{i}(u,X_{u})du\circ dW_{s}^{i} +\sum_{i=1}^{d}\int_{0}^{t}\int_{0}^{s}\partial _{u}\sigma
_{j}(u,X_{u})du\circ
dW_{s}^{i}\\
&+\sum_{i=1}^{d}\int_{0}^{t}\int_{0}^{s}\partial _{\sigma
_{i}}b(u,X_{u})\circ dW_{u}^{i}ds+\int_{0}^{t}\int_{0}^{s}\partial
_{b}b(u,X_{u})duds.
\end{align*}%
We denote
\begin{align*}
&A(t,x)=(\sigma _{i}(t,x),[\sigma _{j},\sigma
_{p}](t,x))_{i,j,p=1,...,d,j\neq p}\quad\mbox{and}\\
&A_{\delta }(t,x)=(\sqrt{\delta }%
\sigma _{i}(t,x),[\sqrt{\delta }\sigma _{j},\sqrt{\delta }\sigma
_{p}](t,x))_{i,j,p=1,...,d,j\neq p}.
\end{align*}
In particular $\lambda _{\ast}(A(t,x))=\lambda (t,x).$

\smallskip

We will need the following estimate for the skeleton $x_t=x_t(\phi)$ as in (\ref{Not3}). And for  $\phi\in L^2([0,T],R^d)$, we set
\begin{equation}\label{eps-phi}
\varepsilon_\phi(\delta)= \Big(\int_0^\delta|\phi_s|^2\,ds\Big)^{1/2}.
\end{equation}

\begin{lemma}\label{lemma-corr}
Let $\delta$ be such that $\varepsilon_\phi(\delta)+\sqrt\delta\leq 1$, $\delta<\frac 1{4n(0,x_0)}$ and
\begin{equation}\label{Corection1}
n(0,x_{0})(\varepsilon
_{\phi }(\delta )+\sqrt{\delta })+\sqrt{\delta }\leq \frac{\sqrt{\lambda(0,x_0)}}{%
8d^{3}n^{2}(0,x_{0})}.
%\varepsilon _{\phi }(\delta )+\sqrt{\delta }\leq \frac{\sqrt{\lambda
%(0,x_{0})}}{8d^{3}n^{3}(0,x_{0})}.
\end{equation}%
Then for every $0\leq t\leq \delta $ and $z\in R^n$,%
\begin{equation}\label{Corection3}
\left\vert z\right\vert^2 _{A_{\delta }(0,x_{0})}\leq 4\left\vert z\right\vert^2
_{A_{\delta }(t,x_{t})}\leq 16\left\vert z\right\vert^2 _{A_{\delta }(0,x_{0})}.
\end{equation}
Moreover,
\begin{equation}
%\sup_{t\leq \delta }\left\vert x_{t}(\phi)-x_{0}-b(0,x_{0})t\right\vert
%_{A_{\delta }(0,x_{0})}\leq 4\varepsilon _{\phi }(\delta )+\sqrt{\delta }.
\sup_{t\leq \delta }\left\vert x_{t}-x_{0}-b(0,x_{0})t\right\vert
_{A_{\delta }(0,x_{0})}\leq 4\varepsilon _{\phi }(\delta )+\frac 1{n(0,x_0)}\delta .
\label{Corection2}
\end{equation}
\end{lemma}

\textbf{Proof}. First, one has $x_s\in B(x_0,1)$ for every $s\leq \delta$. In fact, setting $\tau=\inf\{t>0\,:\,|x_t-x_0|>1\}$, for $s\leq \delta\wedge \tau$ one has
$$
\big|x_s-x_0\big| \leq n(0,x_0)\sqrt\delta(\varepsilon_\phi(\delta)+\sqrt\delta)\leq \frac 12
$$
because $\varepsilon_\phi(\delta)+\sqrt\delta\leq 1$ and $\delta<\frac 1{4n(0,x_0)}$. This gives $s<\tau$. This means that $\delta<\tau$, so that $|x_s-x_0|<1$ for every $s\leq\delta$. Moreover, by using (\ref{Corection1}),
\begin{equation}\label{Corection2bis}
\left\vert x_{s}-x_{0}\right\vert +|s|\leq n(0,x_{0})\sqrt{\delta }(\varepsilon
_{\phi }(\delta )+\sqrt{\delta })+\delta\leq \frac{\sqrt{\lambda(0,x_0))}}{%
8d^{3}n^{2}(0,x_{0})}\times \sqrt{\delta }.
\end{equation}
Now, (\ref{Corection3}) follows immediately from Proposition \ref{NORM3} in Appendix 4 (see page \pageref{page-norm3}).

We prove now (\ref{Corection2}). For $t\leq \delta$, we write now
\begin{equation*}
J_{t}:=x_{t}-x_{0}-b(0,x_{0})t=\int_{0}^{t}(\partial
_{s}x_{s}-b(s,x_{s}))ds+\int_{0}^{t}(b(s,x_{s})-b(0,x_{0}))ds.
%=:J_{t}^{\prime}+J_{t}^{\prime \prime }.
\end{equation*}%
By using inequality (\ref{Norm4bis}) in Lemma \ref{NORM1} from Appendix 4 (see page \pageref{page-norm1}), we get
\begin{align*}
|J_t|^2_{A_\delta(0,x_0)}
&\leq
2t\int_0^t|\partial
_{s}x_{s}-b(s,x_{s})|^2_{A_\delta(0,x_0)}ds
+2t\int_0^t|b(s,x_{s})-b(0,x_0)|^2_{A_\delta(0,x_0)}ds\\
&=:
I'_t
+I''_t
\end{align*}
As for $I'_t$,  we use (\ref{Corection3}): for $s\leq t\leq \delta$ we have
$$
|\partial
_{s}x_{s}-b(s,x_{s})|^2_{A_\delta(0,x_0)}
\leq 4|\partial
_{s}x_{s}-b(s,x_{s})|^2_{A_\delta(s,x_s)}.
$$
Moreover, we can write
$$
\partial _{s}x_{s}-b(s,x_{s}) =\sum_{j=1}^{d}\sigma _{j}(s,x_{s})\phi
_{j}(s)=A_{\delta }(s,x_{s})\psi(s) ,\quad \mbox{with } \psi _{j}(s) =\frac{1}{\sqrt{\delta }}\phi _{j},\  \psi _{i,j}(s)=0
$$
so that
\begin{equation*}
\left\vert \partial _{s}x_{s}-b(s,x_{s})\right\vert _{A_{\delta
}(s,x_{s})}=\left\vert A_{\delta }(s,x_{s})\psi (s)\right\vert _{A_{\delta
}(s,x_{s})}\leq \left\vert \psi (s)\right\vert =\frac{1}{\sqrt{\delta }}%
\left\vert \phi (s)\right\vert .
\end{equation*}%
Then, for $t\leq \delta$ we can write
\begin{equation*}
I'_t\leq
8\delta\int_{0}^{\delta}\left\vert \partial _{s}x_{s}-b(s,x_{s})\right\vert^2 _{A_{\delta
}(s,x_{s})}ds\leq 8\int_{0}^{\delta}\left\vert \phi
(s)\right\vert^2 ds= 8\varepsilon _{\phi }(\delta )^2.
\end{equation*}%
We estimate now $I''_t$: by using (\ref{Corection2bis}),%
\begin{align*}
I''_t&\leq 2\delta\int_0^\delta\frac 1{\lambda_*(A_\delta(0,x_0))}|b(s,x_s)-b(0,x_0)|^2ds\\
&\leq 2\frac{n^2(0,x_0)}{\lambda(0,x_0)}\int_0^t (|s|+|x_s-x_0|)^2ds
\leq
\,\frac{1}{n^{2}(0,x_{0})}\times \delta^2.
\end{align*}
By inserting the estimates for $I'_t$ and $I''_t$, we get
$$
\sup_{t\leq \delta}|J_t|_{A_\delta(0,x_0)}\leq \Big(8\varepsilon_\phi(\delta)^2+\frac 1{n^2(0,x_0)}\delta^2\Big)^{1/2}\leq 4\varepsilon_\phi(\delta)+\frac 1{n(0,x_0)}\delta.
$$
$\square$

\medskip

The main estimate in this section is the following proposition.

\begin{proposition}
\label{DIF1} %Assume that $\delta$ satisfies the assumptions in Lemma \ref{lemma-corr} hold. Moreover,
Let (\ref{Not6}) hold and let $\rho \in (0,1)$ be fixed. Then there exist
some universal constants $C_{1},C_{2}$ (depending on $d$ and $\rho $ only)
such that the following holds. Let $0<\delta \leq R\leq 1$ and $r\in (0,1)$ be
such that
\begin{equation}
\varepsilon _{\phi }(\delta )\leq \frac{r\wedge \sqrt{\lambda (0,x_{0})}}{%
C_{1}n^{3}(0,x_{0})},\quad \delta \leq \frac{r^{5}R}{C_{1}}
\times \frac{\lambda ^{1+3d\rho }(0,x_{0})}{n^{6+6d\rho }(0,x_{0})}
\label{Dif5'}
\end{equation}%
and suppose that%
\begin{equation}
\left\vert X_{0}-x_{0}\right\vert _{A_{\delta }(0,x_{0})}\leq \frac{r}{8}.
\label{Dif7}
\end{equation}%
Then%
\begin{equation}
P\Big(\sup_{t\leq \delta }\left\vert X_{t}-x_{t}\right\vert
_{A_{R}(t,x_{t})}\leq 2r,\left\vert X_{\delta }-x_{\delta }\right\vert
_{A_{\delta }(\delta ,x_{\delta })}\leq r\Big)
\geq \frac{r^{m}}{C_{2}}\exp \Big(-%
\frac{C_{2}n^{2d\rho }(0,x_{0})}{\lambda ^{d\rho }(0,x_{0})}\Big).  \label{Dif7'}
\end{equation}
\end{proposition}

\textbf{Proof}.
For $t\leq \delta$, by using (\ref{Corection3}) we obtain%
\begin{equation*}
\left\vert X_{t}-x_{t}\right\vert _{A_{\delta }(t,x_{t})}\leq 4\left\vert
X_{t}-x_{t}\right\vert _{A_{\delta }(0,x_{0})}\leq 4\sum_{j=1}^{6}\left\vert
I_{j}\right\vert _{A_{\delta }(0,x_{0})}
\end{equation*}%
with%
\begin{eqnarray*}
I_{1} &=&X_{0}-x_{0},\quad I_{2}=Z_{t},\quad I_{3}=R_{t}\quad  \\
I_{4} &=&x_{t}-x_{0}-b(0,x_{0})t,\quad I_{5}=(b(0,X_{0})-b(0,x_{0}))t,
\end{eqnarray*}%
We have to estimate the above terms for $t\leq \delta .$ First%
\begin{equation*}
\left\vert I_{5}\right\vert _{A_{\delta }(0,x_{0})}\leq \frac{n(0,x_{0})}{%
\sqrt{\lambda (0,x_{0})}}\left\vert X_{0}-x_{0}\right\vert \leq \frac{%
n^{2}(0,x_{0})}{\sqrt{\lambda (0,x_{0})}}\left\vert X_{0}-x_{0}\right\vert
_{A_{\delta }(0,x_{0})}\times \sqrt{\delta }\leq \frac{r}{8}
\end{equation*}%
and by (\ref{Corection2})
\begin{equation*}
\left\vert I_{4}\right\vert _{A_{\delta }(0,x_{0})}\leq
4\varepsilon _{\phi }(\delta )+\frac 1{n(0,x_0)}\delta
\leq \frac{r}{8}
\end{equation*}
And by our assumption $\left\vert I_{1}\right\vert _{A_{\delta
}(0,x_{0})}\leq \frac{r}{8}.$ So we have%
\begin{equation*}
\left\vert X_{t}-x_{t}\right\vert _{A_{\delta }(t,x_{t})}\leq \frac{r}{2}%
+4(\left\vert Z_{t}\right\vert _{A_{\delta }(0,x_{0})}+\left\vert
R_{t}\right\vert _{A_{\delta }(0,x_{0})}).
\end{equation*}%
Since $R\geq \delta $, by (\ref{Norm2}) in Lemma \ref{NORM1} from Appendix 4 (see page \pageref{page-norm1}) we have $\left\vert y\right\vert _{A_{R}(0,x_{0})}\leq
\left\vert y\right\vert _{A_{\delta }(0,x_{0})}$ so $\left\vert
I_{i}\right\vert _{A_{R}(0,x_{0})}\leq \frac{r}{8}$ for $i=1,4,5.$ And this
gives%
\begin{equation*}
\left\vert X_{t}-x_{t}\right\vert _{A_{R}(t,x_{t})}\leq \frac{r}{2}%
+4(\left\vert Z_{t}\right\vert _{A_{R}(0,x_{0})}+\left\vert R_{t}\right\vert
_{A_{R}(0,x_{0})}).
\end{equation*}
Using the above inequalities we easily obtain
\begin{align*}
&P\Big(\sup_{t\leq \delta }\left\vert X_{t}-x_{t}\right\vert _{A_{R}(t,x_{t})}
\leq 2r,\left\vert X_{\delta }-x_{\delta }\right\vert _{A_{\delta }(\delta
,x_{\delta })}\leq r\Big) \\
&\quad\geq
P\Big(\sup_{t\leq \delta }\left\vert Z_{t}\right\vert
_{A_{R}(0,x_{0})}+\sup_{t\leq \delta }\left\vert R_{t}\right\vert
_{A_{R}(0,x_{0})} \leq \frac{r}{4},\left\vert Z_{\delta }\right\vert
_{A_{\delta }(0,x_{0})}+\left\vert R_{\delta }\right\vert _{A_{\delta
}(0,x_{0})}\leq \frac{r}{8}\Big)\\
&\quad\geq P\Big(\sup_{t\leq \delta }\left\vert Z_{t}\right\vert _{A_{R}(0,x_{0})} \leq %
\frac{r}{8},\left\vert Z_{\delta }\right\vert _{A_{\delta }(0,x_{0})}\leq
\frac{r}{16}\Big)-P\Big(\sup_{t\leq \delta }\left\vert R_{t}\right\vert _{A_{\delta
}(0,x_{0})}>\frac{r}{8}\Big).
\end{align*}

We upper bound now the last term. First, using the norms inequalities%
\[
P\Big(\sup_{t\leq \delta }\left\vert R_{t}\right\vert _{A_{\delta
}(0,x_{0})}>\frac{r}{8}\Big)\leq P\Big(\sup_{t\leq \delta }\left\vert
R_{t}\right\vert >K\Big)
\]%
with $K=\frac{r\delta \sqrt{\lambda (0,x_{0})}}{8}.$ We define now $\tau
=\inf \{t:\left\vert X_{t}-X_{0}\right\vert \geq \frac{1}{2}).$ Using the
norms inequalties, (\ref{Dif5'}) and (\ref{Dif7}) we obtain $\left\vert
x_{0}-X_{0}\right\vert \leq \frac{1}{2}$ so that for $t\leq \tau $ we have $%
\left\vert X_{t}-x_{0}\right\vert \leq 1.$ It follows that up to $\tau $ the
diffusion process $X$ coincides with a diffusion process $\overline{X}$
which has the coefficients and their derivatives up to order three bounded
by $n(x_{0}).$ We denote by $\overline{R}$ the reminder in which $X$ is
replace with $\overline{X}$ and we write%
\[
P\Big(\sup_{t\leq \delta }\left\vert R_{t}\right\vert >K\Big)\leq P\Big(%
\sup_{t\leq \delta }\left\vert \overline{R}_{t}\right\vert >K\Big)+P(\tau
\leq \delta ).
\]%
Since $\tau =\overline{\tau }:=\inf \{t:\left\vert \overline{X}%
_{t}-X_{0}\right\vert \geq \frac{1}{2})$ a standard reasoning based on
Bernstin's inequality gives $P(\tau \leq \delta )=P(\overline{\tau }\leq
\delta )\leq \exp (-1/C\delta n^{2}(x_{0})).$

In order to estimate the last first we use (\ref{Exp1}) from the Appendix 1
(see Lemma \ref{EXP1} at page \pageref{page-EXP1}) with $k=3,p_{3}=\frac{2}{3%
}$ and with $k=1,p_{1}=2,$ and $K=\frac{r\delta }{4}\sqrt{\lambda (0,x_{0})}.
$ A straightforward computation gives%
\begin{align*}
P\Big(\sup_{t\leq \delta }\left\vert \overline{R}_{t}\right\vert
>\frac{%
r\delta \sqrt{\lambda (0,x_{0})}}{8}\Big)
&\leq
C\exp \Big(-\frac{r^{2}\lambda
(0,x_{0})}{C\delta n^{4}(0,x_{0})}\Big)+C\exp \Big(-\frac{r^{2/3}\lambda
^{1/3}(0,x_{0})}{C\delta ^{1/3}n^{2}(0,x_{0})}\Big) \\
& \leq C\exp \Big(-\frac{r^{2/3}\lambda ^{1/3}(0,x_{0})}{C\delta
^{1/3}n^{2}(0,x_{0})}\Big)
\end{align*}%
the last inequality being a consequence of (\ref{Dif5'}).

Using (\ref{Estimate12})
\begin{equation*}
P\Big(\sup_{t\leq \delta }\left\vert Z_{t}\right\vert _{A_{R}(0,x_{0})}\leq
\frac{r}{8},\left\vert Z_{\delta }\right\vert _{A_{\delta }(0,x_{0})}\leq
\frac{r}{16}\Big)\geq \frac{r^{m}}{2C_{\ast }}\exp \Big(-\frac{C_{\ast }n^{2d\lambda
}(0,x_{0})}{\lambda ^{d\lambda }(0,x_{0})}\Big)
\end{equation*}%
with $C_{\ast }$\ the universal constant in (\ref{Estimate12}). Our
assumption on $\delta $ gives
\begin{equation*}
C\exp \Big(-\frac{r^{2/3}\lambda ^{1/3}(0,x_{0})}{C\delta ^{1/3}n^{2}(0,x_{0})}%
\Big)\leq \frac{1}{2}\times \frac{r^{m}}{2C_{\ast }}\exp \Big(-\frac{C_{\ast
}n^{2d\rho }(0,x_{0})}{\lambda ^{d\lambda }(0,x_{0})}\Big)
\end{equation*}%
so we have proved that%
\begin{equation*}
P\Big(\sup_{t\leq \delta }\left\vert X_{t}-x_{t}\right\vert
_{A_{R}(t,x_{t})}\leq 2r,\left\vert X_{\delta }-x_{\delta }\right\vert
_{A_{\delta }(\delta ,x_{\delta })}\leq r\Big)\geq \frac{r^{m}}{4C_{\ast }}\exp
\Big(-\frac{C_{\ast }n^{2d\rho }(0,x_{0})}{\lambda ^{d\lambda }(0,x_{0})}\Big).
\end{equation*}%
$\square $

\subsection{Chain argument}

We recall that, by the hypothesis (\ref{Not6}) we have some functions $%
\lambda ,n\in L(\mu ,h)$ such that $\lambda (t)\leq 1\wedge \lambda
(t,x_{t}) $ and $n_{t}\geq 1\vee n(t,x_{t})$ such that $\lambda ,n\in L(\mu
,h)$ for some $h>0$ and $\mu \geq 1.$ We also consider some $R,r,\rho \in
(0,1)$ and we define (with $C_{1}$ the constant in (\ref{Dif5'}))
\begin{equation*}
f_{h}(t)=\frac{2}{h}+\frac{C_{1}(\ln \frac{1}{r})^{3}n_{t}^{6+4d\rho }}{%
Rr^{2}\lambda _{t}^{1+d\rho }}+\frac{C_{1}^{2}n_{t}^{6}}{r^{2}\wedge \lambda
_{t}}\left\vert \phi _{t}\right\vert ^{2}.
\end{equation*}%
Notice that, if $d\rho \leq \frac{1}{5}$ then $f_{h}\in L(\mu ^{8},h).$ We
define
\begin{equation*}
\delta (t)=\inf \big\{\delta >0:\int_{t}^{t+\delta }f_{h}(s)ds\geq \frac{1}{\mu
^{8}}\big\}
\end{equation*}

\begin{lemma}
\label{CHAI1}i) One has%
\begin{equation*}
\delta (t)\leq \frac{h}{2}\wedge \frac{Rr_{t}^{2}\lambda _{t}^{1+d\rho }}{%
C_{1}(\ln \frac{1}{r})^{3}n_{t}^{6+4d\rho }},\quad \varepsilon _{\phi
}(\delta (t))\leq \frac{r\wedge \lambda _{t}^{1/2}}{C_{1}n_{t}^{3}}.
\end{equation*}

ii) If $\left\vert t-t^{\prime }\right\vert \leq \delta (t)$ then%
\begin{equation*}
\frac{1}{4\mu ^{16}}\left\vert y\right\vert _{A_{\delta (t)}(t,x_{t})}\leq
\left\vert y\right\vert _{A_{\delta (t^{\prime })}(t^{\prime },x_{t^{\prime
}})}\leq 4\mu ^{8}\left\vert y\right\vert _{A_{\delta (t)}(t,x_{t})}.
\end{equation*}
\end{lemma}

\textbf{Proof}. i) Since $\int_{t}^{t+h/2}\frac{2}{h}ds=1\geq 1/\mu ^{8}$ we
have $\delta (t)\leq \frac{1}{2}h.$ So we may use the properties $L(\mu ,h)$
for $t\leq s\leq t+\delta (t).$ Consequently, for $0<\delta \leq \delta (t)$%
\begin{equation*}
\frac{1}{\mu ^{8}}\geq \int_{t}^{t+\delta }\frac{C_{1}(\ln \frac{1}{r}%
)^{3}n_{s}^{6+4d\rho }}{Rr^{2}\lambda _{s}^{1+d\rho }}ds\geq \frac{1}{\mu
^{8}}\times \frac{C_{1}(\ln \frac{1}{r})^{3}n_{t}^{6+4d\rho }}{Rr^{2}\lambda
_{t}^{1+d\rho }}\times \delta
\end{equation*}%
which gives
\begin{equation*}
\delta (t)\leq \frac{Rr^{2}\lambda _{t}^{1+d\rho }}{C_{1}(\ln \frac{1}{r}%
)^{3}n_{t}^{6+4d\rho }}.
\end{equation*}%
We also have
\begin{equation*}
\frac{1}{\mu ^{8}}\geq \int_{t}^{t+\delta }\frac{C_{1}^{2}n_{s}^{6}}{%
r^{2}\wedge \lambda _{s}}\left\vert \phi _{s}\right\vert ^{2}ds\geq \frac{1}{%
\mu ^{8}}\times \frac{C_{1}^{2}n_{t}^{6}}{r^{2}\wedge \lambda _{t}}%
\int_{t}^{t+\delta }\left\vert \phi _{s}\right\vert ^{2}ds
\end{equation*}%
so that
\begin{equation*}
\varepsilon _{\phi }^{2}(t)\leq \frac{r^{2}\wedge \lambda _{t}}{%
C_{1}^{2}n_{t}^{6}}.
\end{equation*}%
This proves $i).$

$ii)$ We use here next Proposition \ref{NORM3} from Appendix 4 (see page \pageref{page-norm3}).

\smallskip

If $\left\vert t-t^{\prime }\right\vert \leq \delta (t),$ then $%
\left\vert x_{t}-x_{t^{\prime }}\right\vert \leq \delta
^{1/2}(t)(d\varepsilon _{\phi }(\delta (t))+\delta ^{1/2}(t))n_{t}$ so (\ref%
{Norm5}) is verified and we may use (\ref{Norm6}) to obtain
\begin{equation*}
\frac{1}{4}\left\vert y\right\vert _{A_{\delta (t)}(t,x_{t})}\leq \left\vert
y\right\vert _{A_{\delta (t)}(t^{\prime },x_{t^{\prime }})}\leq 4\left\vert
y\right\vert _{A_{\delta (t)}(t,x_{t})}.
\end{equation*}%
It remains to compare $\left\vert y\right\vert _{A_{\delta (t)}(t,x_{t})}$
with $\left\vert y\right\vert _{A_{\delta (t^{\prime })}(t,x_{t})}.$ Since $%
\delta (t^{\prime })\leq \frac{1}{2}h$ and $\left\vert t-t^{\prime
}\right\vert \leq \frac{1}{2}h$ we have $\left\vert t-s\right\vert \leq h$
for every $s\in (t^{\prime },t^{\prime }+\delta (t^{\prime }).$ We use the
property $L(\mu ^{8},h)$ for $f_{h}$ and we obtain%
\begin{equation*}
\mu ^{8}f_{h}(t)\delta (t)\geq \int_{t}^{t+\delta (t)}f_{h}(s)ds=\frac{1}{%
\mu ^{8}}=\int_{t^{\prime }}^{t^{\prime }+\delta (t^{\prime
})}f_{h}(s)ds\geq \mu ^{-8}f_{h}(t)\delta (t^{\prime }).
\end{equation*}%
So $(\delta (t)/\delta (t^{\prime }))^{1/2}\geq \mu ^{-8}.$ Suppose now that
$\delta (t)\leq \delta (t^{\prime }).$ We use then (\ref{Estimate1}) and we
obtain
\begin{equation*}
\frac{1}{\mu ^{16}}\left\vert y\right\vert _{A_{\delta (t)}(t,x_{t})}\leq
\left\vert y\right\vert _{A_{\delta (t^{\prime })}(t,x_{t})}\leq \frac{1}{%
\mu ^{8}}\left\vert y\right\vert _{A_{\delta (t)}(t,x_{t})}.
\end{equation*}%
$\square $

We construct now a time grid in the following way. We put $t_{0}=0$ and%
\begin{equation*}
t_{k}=t_{k-1}+\delta (t_{k-1})
\end{equation*}%
and we denote%
\begin{equation*}
\Theta _{k}=\Big\{\sup_{t_{k-1}\leq t\leq t_{k}}\left\vert
X_{t}-x_{t}\right\vert _{A_{R}(t,x_{t})}\leq r\Big\},\quad \Gamma
_{k}=\Big\{\left\vert X_{t_{k}}-x_{t_{k}}\right\vert _{A_{\delta
(t_{k})}(t_{k},x_{t_{k}})}\leq \frac{r}{8}\Big\}.
\end{equation*}

\begin{proposition}
\label{CHAIN2}i) Suppose that (\ref{Not6}) holds and let $R,r\in (0,1)$ and $%
\rho \in (0,\frac{1}{5d})$. There exists a universal constant $C$
(depending on $d$ and on $\rho )$ such that%
\begin{equation*}
P\big(\cap _{i=1}^{k}\Theta _{i}\cap \Gamma _{i}\big)\geq P\big(\cap _{i=1}^{k-1}\Theta
_{i}\cap \Gamma _{i}\big)\exp \Big(-\frac{Cn_{t_{k-1}}^{2d\rho }}{\lambda
_{t_{k-1}}^{d\rho }}\Big).
\end{equation*}%
ii) Moreover there exists an universal constant $C$ such that%
\begin{align*}
P\big(\sup_{0\leq s\leq T}\left\vert X_{s}-x_{s}\right\vert _{A_{R}(s,x(s))}
\leq r\big)
&\geq \exp \Big(-C\mu ^{9}\int_{0}^{T}f_{h}(t)\frac{n_{s}^{2d\rho }}{%
\lambda _{t}^{d\rho }}dt\Big) \\
&\geq \exp \Big(-C\mu ^{9}\Big(\frac{T}{h}+\frac{1}{r^{2}}\int_{0}^{T}\frac{%
n_{t}^{6+6d\rho }}{\lambda _{t}^{1+2d\rho }}\Big(\frac{(\ln \frac{1}{r})^{3}}{R}%
+\left\vert \phi _{t}\right\vert ^{2}\Big)dt\Big)\Big).
\end{align*}
\end{proposition}

\textbf{Proof }$i)$\textbf{\ }Let%
\begin{equation*}
\widetilde{\Gamma }_{k}=\Big\{\left\vert X_{t_{k}}-x_{t_{k}}\right\vert
_{A_{\delta (t_{k-1})}(t_{k-1},x_{t_{k-1}})}\leq \frac{1}{32\mu ^{8}}r\Big\}.
\end{equation*}%
Using $ii)$ from the previous lemma we obtain $\widetilde{\Gamma }%
_{k}\subset \Gamma _{k}$ so by (\ref{Dif7'})%
\begin{equation*}
P_{t_{k-1}}(\Theta _{k}\cap \Gamma _{k})\geq P_{t_{k-1}}(\Theta _{k}\cap
\widetilde{\Gamma }_{k})\geq \exp \Big(-C\frac{n_{t_{k-1}}^{2d\rho }(x_{0})}{%
\lambda _{t_{k-1}}^{d\rho }(A(x_{0}))}\Big).
\end{equation*}%
The above inequality holds if $\left\vert X_{t_{k-1}}-x_{t_{k-1}}\right\vert
_{A_{\delta (t_{k-1})}(t_{k-1},x_{t_{k-1}})}\leq \frac{r}{8}$ and this is
true on the set $\Gamma _{k-1}.$

$ii)$ Let $N_{T}=\min \{k:t_{k}>T\}.$ Since $X_{0}=x_{0}$ we may use the
recursively the inequality from $i)$ and we obtain%
\begin{equation*}
P\big(\sup_{t\leq T}\left\vert X_{t}-x_{t}\right\vert _{A_{R}(x_{t})}\leq r\big)\geq
P\Big(\cap _{i=1}^{N_{T}}\Theta _{i}\cap \Gamma _{i}\Big)\geq \exp
\Big(-C\sum_{k=1}^{N_{T}}\frac{n_{t_{k-1}}^{2d\rho }}{\lambda _{t_{k-1}}^{d\rho }%
}\Big).
\end{equation*}%
We write%
\begin{eqnarray*}
\int_{0}^{T}f_{h}(s)\frac{n_{s}^{2d\rho }}{\lambda _{s}^{d\rho }}ds &\geq
&\sum_{i=1}^{N_{T}-1}\int_{t_{i-1}}^{t_{i}}f_{h}(s)\frac{n_{s}^{2d\rho }}{%
\lambda _{s}^{d\rho }}ds\geq \frac{1}{\mu ^{3d\rho }}\sum_{i=1}^{N_{T}-1}%
\frac{n_{t_{i-1}}^{2d\rho }}{\lambda _{t_{i-1}}^{d\rho }}%
\int_{t_{i-1}}^{t_{i}}f_{h}(s)ds \\
&=&\frac{1}{\mu ^{8+3d\rho }}\sum_{i=1}^{N_{T}-1}\frac{n_{t_{i-1}}^{2d\rho }%
}{\lambda _{t_{i-1}}^{d\rho }}
\end{eqnarray*}%
the last equality being a consequence of the definition of $\delta (t_{k}).$
$\square $

\section{Appendix 1. Exponential decay for multiple sto\-chastic integrals}\label{app-multiple}

In this section $W=(W^{1},...,W^{d})$ is a standard Brownian motion and $%
\alpha =(\alpha _{1},...,\alpha _{k})\in \{1,...,d\}^{k}$ denotes a multi
index. We use the notation $\overline{\alpha }=(\alpha _{1},...,\alpha
_{k-1}).$\ We consider an adapted and bounded stochastic process $a$ and we
denote by $\left\Vert a\right\Vert _{\infty }$ a constant such that $%
\sup_{t\leq T}\left\vert a(t,\omega )\right\vert \leq \left\Vert a\right\Vert _{\infty }$
almost surely. Then we define the iterated stochastic integrals%
\begin{equation*}
I_{0}(a,W)(t)=a(t),\quad I_{k}^{\alpha }(a,W)(t)=\int_{0}^{t}I_{k-1}^{%
\overline{\alpha }}(a,W)(s)dW_{s}^{\alpha _{k}}.
\end{equation*}

\label{page-EXP1}
\begin{lemma}
\label{EXP1}There exist some universal constants $C_{k},C_{k}^{\prime }$
such that for each $T,K\geq 0$ and every multi-index $\alpha =(\alpha
_{1},...,\alpha _{k})$ one has
\begin{equation}
P\Big(\sup_{t\leq T}\left\vert I_{k}^{\alpha }(a,W)(t)\right\vert \geq K\Big)\leq
C_{k}\exp \Big(-C_{k}^{\prime }\Big( \frac{K}{T^{k/2}\left\Vert a\right\Vert
_{\infty }}\Big) ^{p_{k}}\Big)  \label{Exp1}
\end{equation}%
with%
\begin{equation*}
p_{1}=2,\quad p_{k+1}=\frac{2p_{k}}{2+p_{k}}.
\end{equation*}
\end{lemma}

\textbf{Proof}. We assume that $\left\Vert a\right\Vert _{\infty }=1$ almost
surely (if not we normalize with $\left\Vert a\right\Vert _{\infty })$ and $%
T=1$ (if not we use a scaling argument). We proceed by recurrence. We take
some $Q\geq 0$ and we write%
\begin{eqnarray*}
P(\sup_{t\leq 1}\left\vert I_{k}^{\alpha }(a,W)(t)\right\vert &\geq &K)\leq
I+J\quad with \\
I &=&P(\sup_{t\leq 1}\left\vert I_{k}^{\alpha }(a,W)(t)\right\vert \geq
K,\sup_{t\leq 1}\left\vert I_{k-1}^{\overline{\alpha }}(a,W)(t)\right\vert
\leq Q) \\
J &=&P(\sup_{t\leq 1}\left\vert I_{k-1}^{\overline{\alpha }%
}(a,W)(t)\right\vert \geq Q).
\end{eqnarray*}%
Using the recurrence hypothesis
\begin{equation*}
J\leq C_{k-1}\exp (-C_{k-1}^{\prime }Q^{p_{k-1}}).
\end{equation*}%
We set $h(t)=\int_{0}^{t}\left\vert I_{k-1}^{\overline{\alpha }%
}(a,W)(s)\right\vert ^{2}ds$
and we write $I_{k}^{\alpha }(a,W)(t)=b(h_{t})$ where $b$ is a standard Brownian
motion. So, we obtain%
\begin{equation*}
I\leq P(\sup_{t\leq h_{1}}\left\vert b(t)\right\vert \geq K,h(1)\leq
Q^{2})\leq P(\sup_{t\leq Q^{2}}\left\vert b(t)\right\vert \geq K)\leq C\exp
(-\frac{C^{\prime }K^{2}}{Q^{2}}).
\end{equation*}%
We choose $Q$ solution of $Q^{p_{k-1}}=K^{2}/Q^{2}$ that is $Q=K^{\frac{2}{%
2+p_{k-1}}}.$ Then we obtain
\begin{equation*}
P(\sup_{t\leq 1}\left\vert I_{k}^{\alpha }(a,W)(t)\right\vert \geq K)\leq
C_{k}\exp (-C_{k}^{\prime }K^{\frac{2p_{k-1}}{2+p_{k-1}}})
\end{equation*}%
with $C_{k}=C\vee C_{k-1},C_{k}^{\prime }=C^{\prime }\wedge C_{k-1}^{\prime
}.$ $\square $

\section{Appendix 2. Small perturbations of Gaussian random variables}\label{app-pert}

\subsection{The inverse function theorem}

We give first a quantitative version of the inverse function theorem. We
consider a three time differentiable function
\begin{equation*}
\eta :R^{d}\rightarrow R^{d}\quad \mbox{and}\quad \Phi (\theta ):=\theta +\eta
(\theta ).
\end{equation*}%
We assume that%
\begin{equation*}
\eta (0)=0\quad \mbox{and}\quad \nabla \eta (0)\leq \frac{1}{2}.
\end{equation*}%
In particular this implies that $\nabla \Phi (0)$ is invertible and
\begin{equation*}
\left\vert \nabla \Phi (0)x\right\vert ^{2}\geq \frac{1}{2}\left\vert
x\right\vert ^{2}-\left\vert \nabla \eta (0)x\right\vert ^{2}\geq \frac{1}{2}%
\left\vert x\right\vert ^{2}-\frac{1}{4}\left\vert x\right\vert ^{2}=\frac{1%
}{4}\left\vert x\right\vert ^{2}.
\end{equation*}%
We also have $\left\vert \nabla \Phi (0)x\right\vert \leq \sqrt{3}\left\vert
x\right\vert $ so
\begin{equation*}
\frac{1}{2}\left\vert x\right\vert \leq \left\vert \nabla \Phi
(0)x\right\vert \leq \sqrt{3}\left\vert x\right\vert .
\end{equation*}%
We denote%
\begin{equation*}
c_{2}(\eta )=\max_{i,j=1,d}\sup_{\left\vert x\right\vert \leq 1}\left\vert
\partial _{ij}^{2}\eta (x)\right\vert,\quad c_{3}(\eta
)=\max_{i,j,k=1,d}\sup_{\left\vert x\right\vert \leq 1}\left\vert \partial
_{ijk}^{3}\eta (x)\right\vert
\end{equation*}%
and we take $h_{\eta }>0$ such that
\begin{equation}\label{h-eta}
h_\eta\leq \frac 12 \quad\mbox{and}\quad h_{\eta }\leq \frac{1}{4d^3(c_{2}(\eta )+c_{3}(\eta ))}.
\end{equation}

\begin{proposition}
\label{INVERSE1}
Suppose that $\eta \in C^{3}(R^{d},R^{d})$, $\eta (0)=0$
and $\nabla \eta (0)\leq \frac{1}{2}.$ Then there exists a neighborhood $%
V_{(h_{\eta })}\subset B(0,2h_{\eta })$ of zero such that $\Phi
:V_{(h_{\eta })}\rightarrow B(0,\frac{1}{2}h_{\eta })$ is a  diffeomorphism.
In particular, one has
$$
\Phi^{-1}\ :\ B\Big(0,\frac 12h_{\eta }\Big)\to B(0,2h_{\eta })
$$
and for every $y \in B(0,\frac{1}{2}h_{\eta })$ the following estimates hold:
\begin{equation}  \label{IF1}
\frac{1}{4}\left\vert \Phi ^{-1}(y )\right\vert \leq \left\vert y
\right\vert \leq 4\left\vert \Phi ^{-1}(y )\right\vert.
\end{equation}%
\end{proposition}

\textbf{Proof}. The existence and the differentiability property of the inverse function $\Phi^{-1}$ in a neighborhood of the origin is a well known result from the Inverse Function Theorem. What we aim to prove is that $\Phi^{-1}\ :\ B(0,\frac 12h_{\eta })\to B(0,2h_{\eta })$ and the estimates in (\ref{IF1}).

Since $\eta (0)=0$ we have
\begin{equation*}
\eta (\theta)=\nabla \eta (0)\theta+\int_{0}^{1}(1-t)\left\langle \nabla ^{2}\eta
(t\theta)\theta,\theta\right\rangle dt
\end{equation*}%
with $\nabla ^{2}\eta ^{k}=(\partial _{ij}^{2}\eta
^{k})_{i,j=1,d},k=1,...,d. $ So, given $y\in R^d$ and recalling that
$\nabla \Phi (0)=I+\nabla \eta (0)$, the equation $\Phi
(\theta)=y $ reads%
$$
\theta =U_{y}(\theta),\quad\mbox{with } U_y(\theta):=(\nabla \Phi (0))^{-1}\Big(y
-\int_{0}^{1}(1-t)\left\langle \nabla ^{2}\eta (t\theta)\theta,\theta\right\rangle dt\Big).
$$
Recall that $\frac{1}{2}\left\vert x\right\vert \leq \left\vert \nabla \Phi
(0)x\right\vert .$ Then, for $\theta_1,\theta_2\in B(0,2h_{\eta })$ we have
\begin{eqnarray*}
\left\vert U_{y }(\theta_1)-U_{y}(\theta_2)\right\vert &=&\left\vert (\nabla
\Phi (0))^{-1}\int_{0}^{1}(1-t)(\left\langle \nabla ^{2}\eta
(t\theta_1)\theta_1,\theta_1\right\rangle -\left\langle \nabla ^{2}\eta (t\theta_2)\theta_2,\theta_2\right\rangle
)dt\right\vert \\
&\leq &2\int_{0}^{1}(1-t)\left\vert \left\langle \nabla ^{2}\eta
(t\theta_1)\theta_1,\theta_1\right\rangle -\left\langle \nabla ^{2}\eta (t\theta_2)\theta_2,\theta_2\right\rangle
\right\vert dt \\
&\leq & 2d^3 h_\eta\big (c_2(\eta)+  c_3(\eta)\big)
|\theta_1-\theta_2|\leq \frac 12
\left\vert \theta_1-\theta_2\right\vert,
\end{eqnarray*}%
so that
\begin{equation}\label{IF2bis}
\left\vert U_{y }(\theta_1)-U_{y}(\theta_2)\right\vert
\leq \frac{1}{2}\left\vert \theta_1-\theta_2\right\vert.
\end{equation}
Notice also that for $y \in B(0,\frac{1}{2}h_{\eta })$ and $\theta\in
B(0,2h_{\eta })$\ the above inequality gives%
\begin{equation}
\left\vert U_{y }(\theta)\right\vert \leq \left\vert U_{y
}(\theta)-U_{y }(0)\right\vert +\left\vert U_{y }(0)\right\vert \leq
\frac{1}{2}\left\vert \theta\right\vert +2\left\vert y \right\vert \leq
h_{\eta }+h_{\eta }=2h_{\eta }.  \label{IF3}
\end{equation}%
We define now%
\begin{equation*}
\theta_{0}=0,\quad \theta_{k+1}=U_{y}(\theta_{k}).
\end{equation*}%
From (\ref{IF3}) we know that $\theta_{k}\in B(0,2h_{\eta })$, $k\in N$
and consequently
$$
\left\vert U_{y }(\theta_{k+1})-U_{y
}(\theta_{k})\right\vert \leq \frac{1}{2}\left\vert \theta_{k}-\theta_{k-1}\right\vert .
$$
So the sequence $\theta_{k},k\in N$ converges to the solution of the equation $\theta=U_y(\theta)$, that is $%
\Phi (y)=\theta $. We have thus proved that for any $y \in B(0,\frac{1}{2}
h_{\eta })$ there exists a unique $\theta\in B(0,2h_{\eta })$ such that $%
\Phi (\theta)=y $, that is $\Phi^{-1}\,:\, B(0,\frac 12h_{\eta })\to B(0,2h_{\eta })$ is well defined.

\smallskip

Finally, for $y\in B(0,\frac 12h_{\eta })$ let $\theta=\Phi ^{-1}(y).$ Then $\theta=U_{y}(\theta)$ so, using (\ref{IF3}) $%
\left\vert \theta\right\vert =\left\vert U_{y}(\theta)\right\vert \leq \frac{1}{2}%
\left\vert \theta\right\vert +2\left\vert y \right\vert $ which gives $%
\left\vert \theta\right\vert \leq 4\left\vert y \right\vert .$ Moreover, again by (\ref{IF3}),
\begin{equation*}
\left\vert \theta\right\vert =\left\vert U_{y}(\theta)\right\vert \geq \left\vert
U_{y}(0)\right\vert -\left\vert U_{y}(\theta)-U_{y}(0)\right\vert \geq \frac{1}{2}%
\left\vert y \right\vert -\frac{1}{2}\left\vert \theta\right\vert
\end{equation*}%
which proves that $\left\vert \theta \right\vert \geq \frac 13\left\vert
y\right\vert \geq \frac 14 |y|$. $\square $

\medskip

Let us consider a more specific variant of the local inversion theorem we will need in next Section \ref{app-norms}.
We consider a matrix $B\in \mathcal{M}_{d\times d}$ with columns $B_{i}\in
R^{d},i=1,...,d$ and we suppose that $B$ is invertible. Then we consider the
equation
\begin{equation}
y=B\theta +r(\theta )  \label{Comp1}
\end{equation}%
where $r\in C^{3}(R^{d},R^{d}).$ Our aim is to prove that for small $y$ the
above equation has a unique solution and to obtain some precise estimates
for $\theta $ and its projection on a suitable subspace of $R^d$ in terms of $y.$ In order to do it we have to introduce some
more notations.

We fix $d^{\prime }\in \{1,...,d-1\}$ and we denote $%
d^{\prime \prime }=d-d^{\prime }.$ For $x=(x_{1},...,x_{d})\in R^{d}$ we
denote $\overleftarrow{x}=(x_{1},...,x_{d^{\prime }})$ and $\overrightarrow{x%
}=(x_{d^{\prime }+1},...,x_{d}).$ We denote by $\overleftarrow{B}\in
\mathcal{M}_{d\times d^{\prime }}$ (respectively by $\overrightarrow{B}\in
\mathcal{M}_{d\times d^{\prime \prime }})$ the matrix with columns $%
B_{1},...,B_{d^{\prime }}$ (respectively the matrix with columns $%
B_{d^{\prime }+1},...,B_{d}).$ Let $S=Vect\{B_{1},...,B_{d^{\prime }}).$
Since $B$ is invertible the columns $B_{i},i=1,...,d$ are linearly
independent and so $\dim S=d^{\prime }.$ We denote by $S^{\bot }$ the
orthogonal of $S$ by $\Pi $ the projection on $S$ and by $\Pi ^{\bot }$ the
projection on $S^{\bot }.$ We define $\overrightarrow{B}^{\bot }$ to be the
matrix with columns $B_{i}^{\bot }:=\Pi ^{\bot }B_{i},i=d^{\prime }+1,...,d.$
Since $B_{1},...,B_{d}$ span $R^{d}$ it follows that $B_{i}^{\bot
},i=d^{\prime }+1,...,d$ span $S^{\bot }$ which has dimension $d^{\prime
\prime }.$ So $B_{i}^{\bot },i=d^{\prime }+1,...,d$ are linearly independent.
We conclude that the matrices $B^*B$, $\overleftarrow{B}^{\ast }\overleftarrow{B}\in
\mathcal{M}_{d^{\prime }\times d^{\prime }}$ and $\overrightarrow{B}^{\bot,\ast }\overrightarrow{B}^{\bot}\in
\mathcal{M}_{d^{\prime\prime }\times d^{\prime\prime }}$
are all invertible, and as usual we denote by $\lambda_*(B)$, $\lambda_*(\overleftarrow{B})$ and $\lambda_*(\overrightarrow{B}^\bot)$ the smaller eigenvalue of
$B^*B$, $\overleftarrow{B}^{\ast }\overleftarrow{B}$ and $\overrightarrow{B}^{\bot,\ast }\overrightarrow{B}^{\bot}$ respectively.

\begin{theorem}
\label{INVERSE3}
We assume that the matrix $B$ is invertible and that $%
r(0)=\nabla r(0)=0.$ Set $|B|_\infty=\sup_{i,j=1,\ldots,d}|B_{ij}|$.
Then for every $y\in R^d$ such that
\begin{equation}\label{Comp3}
|y|<\frac {\lambda_*(B)^{1/2}}4\quad\mbox{and}\quad \left\vert y\right\vert < \frac{\lambda_*(B)}{8d^3(c_{2}(r)+c_{2}(r))},
\end{equation}%
the equation (\ref{Comp1}) has a unique solution $\theta $ and%
\begin{equation}
\left\vert \theta \right\vert \leq \frac{4}{\lambda_* ^{1/2}(B)}\left\vert
y\right\vert ,\quad \left\vert \overrightarrow{\theta }\right\vert \leq
\frac{\left\vert B\right\vert _{\infty }}{\lambda_* (\overrightarrow{B}^{\bot
})}\left\vert \Pi ^{\bot }y\right\vert +\frac{16c_{2}(r)\left\vert
B\right\vert _{\infty }}{\lambda_* (\overrightarrow{B}^{\bot })\alpha_* (B)}%
\left\vert y\right\vert ^{2}.  \label{Comp2}
\end{equation}%
In particular if $\left\vert \Pi ^{\bot }y\right\vert \leq \left\vert \Pi
y\right\vert $\ then $\left\vert y\right\vert \leq 2\left\vert \Pi
y\right\vert $ so%
\begin{equation}
\left\vert \theta \right\vert \leq \frac{8}{\lambda_* ^{1/2}(B)}\left\vert \Pi
y\right\vert ,\quad \left\vert \overrightarrow{\theta }\right\vert \leq
\frac{\left\vert B\right\vert _{\infty }}{\lambda_* (\overrightarrow{B}^{\bot
})}\left\vert \Pi ^{\bot }y\right\vert +\frac{64c_{2}(r)\left\vert
B\right\vert _{\infty }}{\lambda_* (\overrightarrow{B}^{\bot })\lambda_* (B)}%
\left\vert \Pi y\right\vert ^{2}.  \label{Comp4}
\end{equation}
\end{theorem}

\textbf{Proof}. We write the equation (\ref{Comp1}) as $B^{-1}y=\theta
+B^{-1}r(\theta )$ and we use Proposition \ref{INVERSE1} with $\eta
(\theta )=B^{-1}r(\theta ).$ Since $\partial _{\alpha }B^{-1}r(\theta
)=B^{-1}\partial _{\alpha }r(\theta )$ we have $c_{2}(\eta )+c_{3}(\eta
)\leq \lambda_* (B)^{-1/2}(c_{2}(r)+c_{3}(r)).$ So our assumption (\ref{Comp3}) ensures
that for some $h_\eta$ fulfilling (\ref{h-eta}), one has $\left\vert B^{-1}y\right\vert \leq \frac{1}{2}h_{\eta }$ and we may
use Proposition \ref{INVERSE1} in order to produce the solution $\theta $ of our
equation. And moreover, by (\ref{IF1}) one has
\begin{equation*}
\left\vert \theta \right\vert \leq 4\left\vert B^{-1}y\right\vert \leq \frac{%
4}{\lambda_* ^{1/2}(B)}\left\vert y\right\vert .
\end{equation*}%
In particular this proves the first inequality in (\ref{Comp2}).\ Using (\ref%
{Comp3}) we also have $\left\vert \theta \right\vert \leq 1.$ Since $%
r(0)=\nabla r(0)=0$ we obtain%
\begin{equation*}
\left\vert r(\theta )\right\vert \leq \max_{\left\vert \alpha \right\vert
=2}\sup_{\left\vert \theta ^{\prime }\right\vert \leq 1}\left\vert \partial
_{\theta }^{\alpha }r(\theta ^{\prime })\right\vert \times \left\vert \theta
\right\vert ^{2}=c_{2}(r)\left\vert \theta \right\vert ^{2}\leq \frac{%
16c_{2}(r)}{\lambda_* (B)}\left\vert y\right\vert ^{2}.
\end{equation*}%
We multiply our equation with $(\overrightarrow{B}^{\bot ,\ast }%
\overrightarrow{B}^{\bot })^{-1}\overrightarrow{B}^{\bot ,\ast }$ and we
obtain%
\begin{equation*}
(\overrightarrow{B}^{\bot ,\ast }\overrightarrow{B}^{\bot })^{-1}%
\overrightarrow{B}^{\bot ,\ast }y=\overrightarrow{\theta }+(\overrightarrow{B%
}^{\bot ,\ast }\overrightarrow{B}^{\bot })^{-1}\overrightarrow{B}^{\bot
,\ast }r(\theta ).
\end{equation*}%
Notice that $\overrightarrow{B}^{\bot ,\ast }y=\overrightarrow{B}^{\bot
,\ast }\Pi ^{\bot }y$ so $\left\vert (\overrightarrow{B}^{\bot ,\ast }%
\overrightarrow{B}^{\bot })^{-1}\overrightarrow{B}^{\bot ,\ast }y\right\vert
\leq \lambda_* ^{-1}(\overrightarrow{B}^{\bot })\left\vert B\right\vert
_{\infty }\left\vert \Pi ^{\bot }y\right\vert $ and this gives%
\begin{equation*}
\left\vert \overrightarrow{\theta }\right\vert \leq \lambda_* ^{-1}(%
\overrightarrow{B}^{\bot })\left\vert B\right\vert _{\infty }\left\vert \Pi
^{\bot }y\right\vert +\lambda_* ^{-1}(\overrightarrow{B}^{\bot })\left\vert
B\right\vert _{\infty }\left\vert r(\theta )\right\vert \leq \frac{%
\left\vert B\right\vert _{\infty }}{\lambda_* (\overrightarrow{B}^{\bot })}%
\left\vert \Pi ^{\bot }y\right\vert +\frac{16c_{2}(r)\left\vert B\right\vert
_{\infty }}{\lambda_* (\overrightarrow{B}^{\bot })\lambda_* (B)}\left\vert
y\right\vert ^{2}.
\end{equation*}%
$\square $

\subsection{Estimates of the density}

For $h>0$ we denote%
\begin{equation}
c_{\ast }(\eta ,h)=\sup_{\left\vert x\right\vert \leq 2h}\max_{i,j}\left\vert \partial _{i}\eta ^{j}(x)\right\vert .  \label{DENSITY2}
\end{equation}%
Let $\Theta $ be a $m$ dimensional centered Gaussian random variable with
covariance matrix $Q.$ We assume that $Q$ is invertible and we denote by $%
\underline{\lambda}(Q)$ and $\overline\lambda(Q)$) the lower and the upper eigenvalue of $Q$ respectively. We also consider a matrix $%
\Gamma \in \mathcal{M}_{n\times m}$ with $n\leq m$ and we recall that $%
\left\vert x\right\vert _{\Gamma }^{2}=\left\langle \Gamma \Gamma ^{\ast
}x,x\right\rangle ,\lambda _{\ast }(\Gamma )$ is the smaller eigenvalue of $%
\Gamma \Gamma ^{\ast }$ and $B_{\Gamma }(y,r)=\{z:\left\vert y-z\right\vert
_{\Gamma }<r\}.$

\begin{lemma}
\label{DENSITY} Suppose that $\Gamma \Gamma ^{\ast }$ is invertible. Let $%
\eta \in C_{b}^{3}(R^{m},R^{n})$ such that $\eta (0)=0$. Set
\begin{equation}
G=y+\Gamma \Theta +\eta (\Theta )  \label{DENSITY1}
\end{equation}%
and assume there exists $r>0$ such that
\begin{equation}
r\leq \frac{1}{2}\lambda _{\ast }(\Gamma )^{1/2}h_{\eta }
%=\frac{\lambda _{\ast }^{1/2}(\Gamma )}{16(c_{2,\eta }+c_{3,\eta })}
\quad \mbox{and}\quad \frac{%
c_{\ast }(\eta ,4r)}{\lambda_*(\Gamma )^{1/2}}\leq \frac{1}{2m},
\label{DENSITY6}
\end{equation}%
$h_\eta$ being defined in (\ref{h-eta}).
Then the law of $G$ has a density $p_{G}$ on $B_{\Gamma }(y,r)$ and for $z\in
B_{\Gamma }(y,r)$ one has%
\begin{align}
\label{DENSITY4}
p_{G}(z)
&\geq
\frac{\underline{\lambda}(Q)^{(m-n)/2}}{\pi^{n/2}\,8^{m/2}\,2^{(m-n)/2}\,(\det Q\,\det \Gamma\Gamma^*)^{1/2}}\,
\exp\Big(-\frac 2{\underline{\lambda}(Q)}\,|z-y|_\Gamma^2\Big)\\
\label{DENSITY5}
p_{G}(z)
&\leq
\frac{\overline{\lambda}(Q)^{(m-n)/2}\,8^{(m-n)/2}\,2^{m/2}}{\pi^{n/2}\,(\det Q\,\det \Gamma\Gamma^*)^{1/2}}\,
\exp\Big(-\frac 1{8\,\overline{\lambda}(Q)}\,|z-y|_\Gamma^2\Big)
\end{align}
In particular,  (\ref{DENSITY4}) and (\ref{DENSITY5}) imply that, for $z\in B_\Gamma(y,r)$,
$$
\Big(\frac{\underline{\lambda}(Q)}{16\,\overline{\lambda}(Q)}\Big)^{m/2}
\,p_{\mathrm{N}(y,\frac 14\underline{\lambda}(Q)\Gamma\Gamma^*)}(z)
\leq p_{G}(z)\leq
\Big(\frac{16\,\overline{\lambda}(Q)}{\underline{\lambda}(Q)}\Big)^{m/2}
\,p_{\mathrm{N}(y,4\overline{\lambda}(Q)\Gamma\Gamma^*)}(z)
$$
where $p_{\mathrm{N}(y,BB^*)}$ denotes the Gaussian density with mean $y$ and covariance matrix $BB^*$.

\end{lemma}

\textbf{Proof. Step 1.} We assume first that $n=m,y=0$ and $\Gamma $ is the
identity matrix. We denote $\Phi (\theta )=\theta +\eta (\theta )$, so that $\Phi(\Theta)=G.$ Let $%
f:R^{m}\rightarrow R$ be a non negative measurable function with the support
included in the (Euclidian) ball $B(0, r)$, with $r$ fulfilling (\ref{DENSITY6}). Using a
change of variable and Proposition \ref{INVERSE1}, we obtain
\begin{align*}
E(f(\Phi (\Theta ))) &=\int_{\{\theta\in \Phi ^{-1}(B(0,r))\}}f(\Phi (\theta))\frac{1}{(2\pi )^{m/2}(\det Q)^{1/2}}\exp (-\frac{1}{2}%
\left\langle Q^{-1}\theta,\theta\right\rangle )d\theta \\
&=\int_{B(0, r)}f(z)\,p_{\Phi(\Theta)}(z)
\,dz,
\end{align*}%
where we have set, for $z\in B(0,r)$,
$$
p_{\Phi(\Theta)}(z)
=\frac{1}{(2\pi )^{m/2}\left\vert
\det \nabla \Phi (\Phi ^{-1}(z))\right\vert (\det Q)^{1/2}}\,
\exp (-\frac{1}{2}%
\left\langle Q^{-1}\Phi ^{-1}(z),\Phi ^{-1}(z)\right\rangle )
$$
Since $r\leq h_\eta$, if $z\in B(0,r)$ one has $\theta=\Phi ^{-1}(z)\in B(0,4r)$ and for $x\in B(0,4r)$ we have%
\begin{equation*}
\frac{1}{2}\left\vert x\right\vert ^{2}\leq (1-mc_{\ast }(\eta ,h_{\eta
}))\left\vert x\right\vert ^{2}\leq \left\vert \left\langle \nabla \Phi
(\theta)x,x\right\rangle \right\vert \leq (1+mc_{\ast }(\eta ,h_{\eta
}))\left\vert x\right\vert ^{2}\leq 2\left\vert x\right\vert ^{2},
\end{equation*}%
because $c_{\ast }(\eta ,4r)\leq \frac 1{2m}$.
Therefore, if $z\in B(0,r)$ then
$$2^{-m}\leq
\left\vert \det \nabla \Phi (\Phi ^{-1}(z))\right\vert \leq 2^{m}.
$$
Moreover, using (\ref{IF1}) we obtain%
\begin{eqnarray*}
\left\langle Q^{-1}\Phi ^{-1}(z),\Phi ^{-1}(z)\right\rangle &\leq &\frac{1}{%
\underline{\lambda}(Q)}\left\vert \Phi ^{-1}(z)\right\vert ^{2}\leq \frac{4}{%
\underline{\lambda}(Q)}\left\vert z\right\vert ^{2}\quad \mbox{and} \\
\left\langle Q^{-1}\Phi ^{-1}(z),\Phi ^{-1}(z)\right\rangle &\geq &\frac{1}{%
\overline{\lambda}(Q)}\left\vert \Phi ^{-1}(z)\right\vert ^{2}\geq \frac{1}{%
4\overline{\lambda}(Q)}\left\vert z\right\vert ^{2}
\end{eqnarray*}%
So, as $z\in B(0,r)$ we get
\begin{equation}\label{DENSITY3}
\frac{1}{(8\pi )^{m/2}\sqrt{\det Q}}\exp (-\frac{2}{\underline{\lambda}(Q)}%
\left\vert z\right\vert ^{2})\leq p_{\Phi (\Theta )}(z)\leq \frac{2^{m/2}}{%
\pi ^{m/2}\sqrt{\det Q}}\exp (-\frac{1}{8\overline{\lambda} (Q)}\left\vert
z\right\vert ^{2})
\end{equation}

\textbf{Step 2}. We still assume that $n=m$ but now $y$ and $\Gamma $ are
general, with $\Gamma$ invertible. We write $G=y+\Gamma (\theta + \eta_\Gamma(\theta ))$ with $%
\eta_\Gamma(\theta )=\Gamma ^{-1}\eta (\theta )$ and denote $\Phi_\Gamma
(\theta )=\theta + \eta_\Gamma(\theta ).$ One has
$c_2(\eta_\Gamma)+c_3(\eta_\Gamma)\leq \lambda_*(\Gamma)^{-1/2}(c_2(\eta)+c_3(\eta))$, so that  $h_{\eta_\Gamma}\geq \lambda _{\ast }(\Gamma )^{1/2}h_{\eta }$ and then (\ref{DENSITY6}) gives $r\leq \frac 12 h_{\eta_\Gamma}$. Moreover, since $c_*(\eta_\Gamma, 4r)\leq \lambda_*(\Gamma)^{-1/2}c_*(\eta, 4r)$, (\ref{DENSITY6}) gives also
$c_*(\eta_\Gamma, 4r)\leq \frac 1{2m}$.
And since $\left\vert \Gamma x\right\vert _{\Gamma }=\left\vert
x\right\vert $, one has
$G\in B_{\Gamma }(y,r)$
iff $\Phi_\Gamma (\Theta )\in B(0,r)$. Then
by a change of variable, for $z\in B_{\Gamma }(y,r)$ we
have%
\begin{equation*}
p_{G}(z)=\frac{1}{\left\vert \det \Gamma \right\vert }p_{\Phi_\Gamma (\Theta
)}(\Gamma ^{-1}(z-y)).
\end{equation*}%
Since $\left\vert \Gamma ^{-1}(z-y)\right\vert =\left\vert z-y\right\vert
_{\Gamma }$ we use (\ref{DENSITY3}) and we obtain%
\begin{align*}
p_{G}(z)\geq &\frac{1}{(8\pi )^{m/2}\sqrt{\det Q}\left\vert \det \Gamma \right\vert }%
\exp (-\frac{2}{\underline{\lambda}(Q)}\left\vert z-y\right\vert _{\Gamma
}^{2}) \\
p_{G}(z)\leq &\frac{2^{m/2}}{\pi ^{m/2}\sqrt{\det Q}\left\vert \det
\Gamma \right\vert }\exp (-\frac{1}{8\overline{\lambda}(Q)}\left\vert
z-y\right\vert _{\Gamma }^{2}).
\end{align*}

\textbf{Step 3}. Now we allow $n$ to be strictly smaller than $m.$ Since $%
\Gamma \Gamma ^{\ast }$ is invertible the lines $\Gamma ^{1},...,\Gamma
^{n}\in R^{m}$ of $\Gamma $ are linearly independent. We denote $%
S=Vect\{\Gamma ^{1},...,\Gamma ^{n}\}$ and we take $\Gamma ^{n+1},...,\Gamma
^{m}$ to be an orthonormal basis in the orthogonal of $S.$ Then we define $%
\widetilde{\Gamma }\in \mathcal{M}_{m\times m}$ to be the matrix with lines $%
\Gamma ^{1},...,\Gamma ^{n},\Gamma ^{n+1},...,\Gamma ^{m}$ and we notice
that
\begin{equation*}
\widetilde{\Gamma }\widetilde{\Gamma }^{\ast }=\left(
\begin{tabular}{ll}
$\Gamma \Gamma ^{\ast }$ & $0$ \\
$0$ & $I_{m-n}$%
\end{tabular}%
\right)
\end{equation*}%
where $I_{m-n}\in \mathcal{M}_{(m-n)\times (m-n)}$ is the identity matrix$.$
In particular %$\lambda _{\ast }(\widetilde{\Gamma })=\lambda _{\ast }(\Gamma )\wedge 1$,
$\det \widetilde{\Gamma }=%
\sqrt{\det \Gamma \Gamma ^{\ast }}$ and for $z=(z_{1},z_{2}),z_{1}\in
R^{n}, $ $z_{2}\in R^{m-n}$ we have $\left\vert z\right\vert _{\widetilde{%
\Gamma }}^{2}=\left\vert z_{1}\right\vert _{\Gamma }^{2}+\left\vert
z_{2}\right\vert ^{2}.$ Moreover, for $y\in R^{n}$ we denote $\widetilde{y}%
=(y,0)$ and we also set $\widetilde\eta(\theta)=(\eta(\theta),0)$. So, we define $H=\widetilde{y}+\widetilde{\Gamma }\Theta +\widetilde{%
\eta }(\Theta )$, and we notice that $h_{\widetilde\eta}=h_\eta$ and $c_*(\widetilde\eta,4r)=c_*(\eta,4r)$. For the density of $H=(H_1,H_2)$ we can use the estimate from the previous step. Notice that since $H_2$ is an orthogonal transformation of a Gaussian random variable, one easily gets that the estimates hold for $(z,u)\in R^m$ such that $z\in B(0,r)$ and $u\in R^{m-n}$.
Now, since $H_{1}=G$ we obtain
\begin{align*}
p_{G}(z) =&\int_{R^{m-n}}p_{H}(z,u)du\\
\geq &\int_{R^{m-n}}\frac{1}{(8\pi )^{d/2}\sqrt{\det Q%
}\left\vert \det \widetilde{\Gamma }\right\vert }\exp (-\frac{2}{\underline{\lambda}
(Q)}\Big(\left\vert z-y\right\vert _{\Gamma }^{2}+\left\vert
u\right\vert ^{2}\Big)du \\
=&\frac{(4\lambda _{\ast }(Q))^{\frac{1}{2}(m-n)}}{(8\pi )^{m/2}\sqrt{\det Q%
}\sqrt{\left\vert \det \Gamma \Gamma ^{\ast }\right\vert }}\exp \Big(-\frac{2}{%
\underline{\lambda} (Q)}\left\vert z-y\right\vert _{\Gamma }^{2}\Big).
\end{align*}%

The proof of the other inequality is the same. $\square $

\section{Appendix 3. Support Property}\label{app-supp}

In this section we prove (\ref{Estimate5'}). Let $B=(B^{1},...,B^{d-1})$ be
a standard Brownian motion. We consider the analogues of the covariance
matrix $Q_{i}(B)$ considered in the previous sections: we define a symmetric
square matrix of dimension $d\times d$ by
\begin{eqnarray*}
Q^{d,d} &=&1,\quad Q^{d,j}=Q^{j,d}=\int_{0}^{1}B_{s}^{j}ds,\quad j=1,...,d-1,
\\
Q^{j,p} &=&Q^{p,j}=\int_{0}^{1}B_{s}^{j}B_{s}^{p}ds,\quad j,p=1,...,d-1
\end{eqnarray*}%
and we denote by $\underline{\lambda}(Q)$ (respectively by $\overline{\lambda}(Q)) $ the lower (respectively larger) eigenvalue of $Q$.

%\begin{remark}
%In order to avoid notational confusion we recall that for a general matrix $%
%M $ we have denoted $\lambda _{\ast }(M)$ the lower eigenvalue of $MM^{\ast
%}.$ With this notation we have $\lambda _{\ast }(B)=\lambda _{\ast
%}^{1/2}(Q).$ But this does not come on in our reasoning.
%\end{remark}

For a measurable function $g:[0,1]\rightarrow R^{d-1}$ we denote%
\begin{eqnarray*}
\alpha_{g}(\xi ) &=&\xi _{d}+\int_{0}^{1}\left\langle g_{s},\xi _{\ast
}\right\rangle ds,\quad \beta_{g}(\xi )=\int_{0}^{1}\left\langle g_{s},\xi
_{\ast }\right\rangle ^{2}ds-\left( \int_{0}^{1}\left\langle g_{s},\xi
_{\ast }\right\rangle ds\right) ^{2}\quad with \\
\xi &=&(\xi _{1},...,\xi _{d})\in R^{d}\quad and\quad \xi _{\ast }=(\xi
_{1},...,\xi _{d-1}).
\end{eqnarray*}

We need the following two preliminary lemmas.

\begin{lemma}
\label{SUPORT1}With $g(s)=B_{s},s\in \lbrack 0,1]$ we have%
\begin{equation*}
\left\langle Q\xi ,\xi \right\rangle =\alpha_{B}^{2}(\xi )+\beta_{B}(\xi ).
\end{equation*}%
As a consequence, one has
\begin{equation*}
\underline{\lambda}(Q)=\inf_{\left\vert \xi \right\vert =1}(\alpha_{B}^{2}(\xi
)+\beta_{B}(\xi ))\quad and\quad \overline{\lambda}(Q)\leq \sup_{\left\vert \xi
\right\vert =1}(\alpha_{B}^{2}(\xi )+\beta_{B}(\xi ))\leq
\big(1+\sup_{t\leq 1}\left\vert
B_{t}\right\vert\big)^2.
\end{equation*}%
Taking $\xi _{\ast }=0$ and $\xi _{d}=1$ we obtain $\left\langle Q\xi ,\xi
\right\rangle =1$ so that $\underline{\lambda}(Q)\leq 1\leq \overline{\lambda}(Q). $
\end{lemma}

\textbf{Proof}. By direct computation%
\begin{eqnarray*}
\left\langle Q\xi ,\xi \right\rangle &=&\xi _{d}^{2}+2\xi
_{d}\int_{0}^{1}\left\langle B_{s},\xi _{\ast }\right\rangle ds+\left(
\int_{0}^{1}\left\langle B_{s},\xi _{\ast }\right\rangle ds)\right) ^{2} \\
&&+\int_{0}^{1}\left\langle B_{s},\xi _{\ast }\right\rangle ^{2}ds-\left(
\int_{0}^{1}\left\langle B_{s},\xi _{\ast }\right\rangle ds\right) ^{2} \\
&=&\left( \xi _{d}+\int_{0}^{1}\left\langle B_{s},\xi _{\ast }\right\rangle
ds\right) ^{2}+\int_{0}^{1}\left\langle B_{s},\xi _{\ast }\right\rangle
^{2}ds-\left( \int_{0}^{1}\left\langle B_{s},\xi _{\ast }\right\rangle
ds\right) ^{2}.
\end{eqnarray*}%
The remaining statements follow straightforwardly.
$\square $

\begin{proposition}
\label{SUPORT2}For each $p\geq 1$ one has%
\begin{equation}
E(\left\vert \det Q\right\vert ^{-p})\leq C_{p,d}<\infty  \label{Sup1}
\end{equation}%
where $C_{p,d}$ is a constant depending on $p,d$ only.
\end{proposition}

\textbf{Proof}. By Lemma 7-29, pg 92 in \cite{[BGJ]}, for every $p\in (0,\infty )$
one has%
\begin{equation*}
\frac{1}{\left\vert \det Q\right\vert ^{p}}\leq \frac{1}{\Gamma (p)}%
\int_{R^{d}}\left\vert \xi \right\vert ^{d(2p-1)}e^{-\left\langle Q\xi ,\xi
\right\rangle }d\xi .
\end{equation*}%
Let $\theta (\xi _{\ast }):=\int_{0}^{1}\left\langle B_{s},\xi _{\ast
}\right\rangle ds.$\ Using the previous lemma%
\begin{eqnarray*}
\int_{R^{d}}\left\vert \xi \right\vert ^{d(2p-1)}e^{-\left\langle Q\xi ,\xi
\right\rangle }d\xi &=&\int_{R^{d}}(\xi _{d}^{2}+\left\vert \xi _{\ast
}\right\vert ^{2})^{d(2p-1)/2}e^{-(\xi _{d}+\theta (\xi _{\ast
}))^{2}-\beta_{B}(\xi _{\ast })}d\xi \\
&\leq &C\int_{R^{d-1}}((1+\theta ^{2}(\xi _{\ast }))^{d(2p-1)/2}+\left\vert
\xi _{\ast }\right\vert ^{d(2p-1)})e^{-\beta_{B}(\xi _{\ast })}d\xi _{\ast } \\
&\leq &C\int_{R^{d-1}}\sup_{t\leq 1}1\vee \left\vert B_{t}\right\vert
^{d(2p-1)}(1+\left\vert \xi _{\ast }\right\vert ^{d(2p-1)+1})e^{-\beta_{B}(\xi
_{\ast })}d\xi _{\ast }.
\end{eqnarray*}%
We integrate and we use Schwartz inequality in order to obtain%
\begin{equation*}
E\Big(\frac{1}{\left\vert \det Q\right\vert ^{p}}\Big)\leq C+C\int_{\{\left\vert \xi
_{\ast }\right\vert \geq 1\}}(E((1+\left\vert \xi _{\ast }\right\vert
^{d(2p-1)+1})^{2}e^{-2\beta_{B}(\xi _{\ast })}))^{1/2}d\xi _{\ast }.
\end{equation*}%
For each fixed $\xi _{\ast }$ the process $b_{\xi _{\ast }}(t):=\left\vert
\xi _{\ast }\right\vert ^{-1}\left\langle B_{t},\xi _{\ast }\right\rangle $
is a standard Brownian motion and $\beta_{B}(\xi _{\ast })=\left\vert \xi _{\ast
}\right\vert ^{2}\int_{0}^{1}(b_{\xi _{\ast }}(t)-\int_{0}^{1}b_{\xi _{\ast
}}(s)ds)^{2}dt=:\left\vert \xi _{\ast }\right\vert ^{2}V_{\xi _{\ast }}$\
where $V_{\xi _{\ast }}$\ is the variance of $b_{\xi _{\ast }}$\ with
respect to the time.\ Then it is proved in \cite{[DY]} (see (1.f), p. 183) that
\begin{equation*}
E(e^{-2\beta_{B}(\xi _{\ast })})=E(e^{-2\left\vert \xi _{\ast }\right\vert
^{2}V_{\xi _{\ast }}})=\frac{2\left\vert \xi _{\ast }\right\vert ^{2}}{\sinh
2\left\vert \xi _{\ast }\right\vert ^{2}}.
\end{equation*}%
We insert this in the previous inequality and we obtain $E(\left\vert \det
Q\right\vert ^{-p})<\infty .$ $\square $

\medskip

We are now able to give the main result in this section. We define%
\begin{equation}
q(B)=\sum_{i=1}^{d-1}\left\vert B_{1}^{i}\right\vert +\sum_{j\neq
p}\left\vert \int_{0}^{1}B_{s}^{j}dB_{s}^{p}\right\vert  \label{Sup4}
\end{equation}%
and for $\varepsilon ,\rho >0$ de denote
\begin{equation}
\Lambda _{\rho ,\varepsilon }(B)=\{\det Q(B)\geq \varepsilon ^{\rho
},\sup_{t\leq 1}\left\vert B_{t}\right\vert \leq \varepsilon ^{-\rho
},q(B)\leq \varepsilon \}.  \label{Sup3}
\end{equation}

\begin{proposition}
\label{SUPORT3}There exist some universal constants $c_{\rho
,d},\varepsilon _{\rho ,d}\in (0,1)$ (depending on $\rho $ and $d$ only)
such that for every $\varepsilon \in (0,\varepsilon _{\rho ,d})$ one has%
\begin{equation}
P(\Lambda _{\rho ,\varepsilon }(B))\geq c_{\rho ,d}\times \varepsilon ^{%
\frac{1}{2}d(d+1)}.  \label{Sup2}
\end{equation}
\end{proposition}

\textbf{Proof.} Using the previous proposition and Chebyshev's inequality we get
$$
P(\det Q<\varepsilon ^{\rho })\leq \varepsilon ^{p\rho }E\left\vert \det
Q\right\vert ^{-p}\leq C_{p,d}\varepsilon ^{p\rho }
\quad\mbox{and}\quad%
P(\sup_{t\leq 1}\left\vert B_{t}\right\vert >\varepsilon ^{-\rho })\leq \exp
(-\frac{1}{C\varepsilon ^{2\rho }}).
$$
Let $q^{\prime }(B)=\sum_{i=1}^{d-1}\left\vert B_{1}^{i}\right\vert
+\sum_{j<p}\left\vert \int_{0}^{1}B_{s}^{j}dB_{s}^{p}\right\vert .$ Since $%
\left\vert \int_{0}^{1}B_{s}^{j}dB_{s}^{p}\right\vert \leq \left\vert
B_{1}^{j}\right\vert \left\vert B_{1}^{p}\right\vert +\left\vert
\int_{0}^{1}B_{s}^{p}dB_{s}^{j}\right\vert $ we have $q(B)\leq 2q^{\prime
}(B)+q'(B)^2$ so that $\{q^{\prime }(B)\leq \frac{1}{3}\varepsilon \}\subset
\{q(B)\leq \varepsilon \}.$ We will now use the following fact: consider the
diffusion process $X=(X^{i},X^{j,p},$ $i=1,...,d,1\leq j<p\leq d)$ solution
of the equation $dX_{t}^{i}=dB_{t}^{i},dX_{t}^{j,p}=X_{t}^{j}dB_{t}^{p}.$
The strong H\"{o}rmander condition holds for this process and the support of
the law of $X_{1}$ is the whole space. So the law of $X_{1}$ is absolutely
continuous with respect to the Lebesgue measure and has a continuous and
strictly positive density $p.$ This result is well known (see for example
\cite{[KS]} or \cite{[BC]}). We denote $c_{d}:=\inf_{\left\vert x\right\vert \leq 1}p(x)>0$
and this is a constant which depends on $d$ only. Then, by observing that $q'(B)\leq \sqrt{m} \,|X_1|$, where $m=\frac{1}{2}d(d+1)$ is the dimension of the diffusion $X$, we get \begin{equation*}
P(q(B)\leq \varepsilon )\geq P\Big(q^{\prime }(B)\leq \frac{\varepsilon }{3}%
\Big)\geq P\Big(\left\vert X_{1}\right\vert \leq \frac{\varepsilon }{3\sqrt m}\Big)\geq \frac{%
\varepsilon ^{m}}{(3\sqrt m)^{m}}\times \bar c_{d},
\end{equation*}%
with $\bar c_d>0$. So finally we obtain
\begin{equation*}
P(\Lambda _{\rho ,\varepsilon }(B))\geq \bar c_{d}\varepsilon ^{\frac{1}{2}%
d(d+1)}-C_{p,d}\varepsilon ^{p\rho }-\exp (-\frac{1}{C\varepsilon ^{2\rho }}%
).
\end{equation*}%
Choosing $p>\frac{1}{2\rho }d(d+1)$ and $\varepsilon $ small we obtain our
inequality$.\square $

\section{Appendix 4. Norms and distances}\label{app-norms}

In this section we use the notation from Section 3 and 4. We consider the
matrix $A$ with columns $a_{i},[a]_{j,p}=a_{j,p}-a_{p,j},i=1,...,d,j\neq p$
defined in Section 3.2 and we assume that the non degeneracy condition (\ref%
{Estimate0}) holds. For notational convenience we denote $%
A_{i}=a_{i},i=1,...,d$ and $A_{i},i=d+1,...,m$ will be an enumeration of $%
[a]_{j,p},1\leq j,p\leq d,j\neq p.$ We work with the norm $\left\vert
y\right\vert _{A_{R}}^{2}=\left\langle (A_{R}A^*_{R})^{-1}y,y\right\rangle$, $y\in R^n$.
We have the following simple properties:

\label{page-norm1}
\begin{lemma}
\label{NORM1}

\begin{itemize}
\item[$i)$] For every $y\in R^n$ and $0<R\leq R^{\prime }\leq 1$ one has
\begin{align}
\sqrt{\frac{R}{R^{\prime }}}\left\vert y\right\vert _{A_{R}}&\geq \left\vert
y\right\vert _{A_{R^{\prime }}}\geq \frac{R}{R^{\prime }}\left\vert
y\right\vert _{A_{R}}\quad\mbox{and}  \label{Norm2}\\
\frac{1}{\sqrt{R}\sqrt{\lambda ^{\ast }(A)}}\left\vert y\right\vert &\leq
\left\vert y\right\vert _{A_{R}}\leq \frac{1}{R\sqrt{\lambda _{\ast }(A)}}%
\left\vert y\right\vert .  \label{Norm3}
\end{align}%
\item[$ii)$] For every $z\in R^{m}$ and $R>0$ one has
\begin{equation}
\left\vert A_{R}z\right\vert _{A_{R}}\leq \left\vert z\right\vert .
\label{Norm4}
\end{equation}
\item[$iii)$] For every $\mu\in L^2([0,T];R^m)$ and $R>0$ one has
\begin{equation}
\Big|\int_0^t\mu_s\,ds\Big|^2_{A_R}
\leq t\int_0^t|\mu_s|^2_{A_R}\,ds,\quad t\in [0,T].
\label{Norm4bis}
\end{equation}
\end{itemize}
\end{lemma}

\textbf{Proof}. $i)$ It is easy to check that
\begin{equation*}
\frac{R^{\prime }}{R}A_{R}A_{R}^{\ast }\leq A_{R^{\prime }}A_{R^{\prime
}}^{\ast }\leq \left( \frac{R^{\prime }}{R}\right) ^{2}A_{R}A_{R}^{\ast }
\end{equation*}%
which is equivalent with (\ref{Norm2}). This also implies (one takes $%
R^{\prime }=1$ so $A_{R^{\prime }}=A)$ that%
\begin{equation*}
\frac{1}{R}\lambda _{\ast }(A_{R})\leq \lambda _{\ast }(A)\leq \frac{1}{R^{2}%
}\lambda _{\ast }(A_{R})\quad and\quad \frac{1}{R}\lambda ^{\ast
}(A_{R})\leq \lambda ^{\ast }(A)\leq \frac{1}{R^{2}}\lambda ^{\ast }(A_{R})
\end{equation*}%
which immediately gives (\ref{Norm3}).

\smallskip

$ii)$ For $z\in R^m$, we write $z=A_{R}^{\ast }y+w$ with $y\in
R^{n}$ and $w\in (\mathrm{Im}A_{R}^{\ast })^{\bot }=KerA_{R}.$ Then $%
A_{R}z=A_{R}A_{R}^{\ast }y$ so that%
\begin{eqnarray*}
\left\vert A_{R}z\right\vert _{A_{R}}^{2} &=&\left\vert A_{R}A_{R}^{\ast
}y\right\vert _{A_{R}}^{2}=\left\langle (A_{R}A_{R}^{\ast
})^{-1}A_{R}A_{R}^{\ast }y,A_{R}A_{R}^{\ast }y\right\rangle \\
&=&\left\langle z,A_{R}A_{R}^{\ast }z\right\rangle =\left\langle A_{R}^{\ast
}y,A_{R}^{\ast }y\right\rangle =\left\vert A_{R}^{\ast }y\right\vert
^{2}\leq \left\vert z\right\vert ^{2}
\end{eqnarray*}
and (\ref{Norm4}) holds.

\smallskip

$iii)$ For $\mu\in L^2([0,T];R^m)$ and $t\in [0,T]$ one has
\def\<{\langle}
\def\>{\rangle}
\begin{align*}
\Big|\int_0^t\mu_sds\Big|^2_{A_R}
=&\big\<{A_R}^{-1}\int_0^t\mu_sds,\int_0^t\mu_sds\big\>
=\int_0^t\int_0^t\big\<{A_R}^{-1}\mu_s,\mu_u\big\>dsdu\\
=&\frac 12\int_0^t\int_0^t\Big(\big\<{A_R}^{-1}(\mu_s-\mu_u),\mu_s-\mu_u\big\>
-\big\<{A_R}^{-1}\mu_s,\mu_s\big\>
-\big\<{A_R}^{-1}\mu_u,\mu_u\big\>\Big)dsdu\\
%&\qquad\quad-\big\<{A_R}^{-1}\mu_s,\mu_s\big\>
%-\big\<{A_R}^{-1}\mu_u,\mu_u\big\>\Big)dsdu\\
=&\frac 12\int_0^t\int_0^t\Big(|\mu_s-\mu_u|^2_{{A_R}}-2|\mu_s|^2_{{A_R}}\Big)dsdu\\
\leq & \int_0^t\int_0^t|\mu_u|^2_{A_R}dsdu
=  t\int_0^t|\mu_u|^2_{A_R}du.
\end{align*}
$\square $

\medskip

We give now some lower and upper bounds for $\left\vert y\right\vert
_{A_{R}}.$ We denote $S=Vect\{A_{i},i=1,...,d\}$ and $\Pi _{S}$ is the
projection on $S.$ $S^{\bot }$ is the orthogonal of $S$ and $\Pi _{S^{\bot
}} $ is the projection on $S^{\bot }.$ Moreover we denote
\begin{eqnarray}
\underline{\lambda} _{S} &=&\inf_{\xi \in S,\left\vert \xi \right\vert
=1}\sum_{i=1}^{d}\left\langle A_{i},\xi \right\rangle ^{2},\quad \overline{\lambda}
_{S}%^{\ast }
=\sup_{\xi \in S,\left\vert \xi \right\vert
=1}\sum_{i=1}^{d}\left\langle A_{i},\xi \right\rangle ^{2}  \label{Norm1a} \\
\underline{\lambda} _{S^{\bot }} &=&\inf_{\xi \in S^{\bot },\left\vert \xi \right\vert
=1}\sum_{i=d+1}^{m}\left\langle A_{i},\xi \right\rangle ^{2},\quad \overline{\lambda}
_{S^{\bot }}=\sup_{\xi \in S^{\bot },\left\vert \xi \right\vert
=1}\sum_{i=d+1}^{m}\left\langle A_{i},\xi \right\rangle ^{2}.  \notag
\end{eqnarray}%
By the very definition $\underline{\lambda} _{S}>0$ and under assumption (\ref{Estimate0}%
) we also have $\underline{\lambda} _{S^{\bot }}>0.$ And $\overline{\lambda} _{S}\leq
\lambda ^{\ast }(A),\overline{\lambda} _{S^{\bot }}\leq \lambda ^{\ast }(A).$

\begin{proposition}
\label{NORM2}Suppose that (\ref{Estimate0}) holds and let
\begin{equation}
R< \frac{\underline{\lambda} _{S}}{4\sum_{i=d+1}^{m}\left\vert \Pi
_{S}A_{i}\right\vert ^{2}}.  \label{Norm2a}
\end{equation}%
Then for every $y\in R^{n}$%
\begin{equation}
\frac{1}{4R\,\overline{\lambda }_{S}}\left\vert \Pi _{S}y\right\vert ^{2}+\frac{1%
}{4R^{2}\,\overline{\lambda} _{S^{\bot }}}\left\vert \Pi _{S^{\bot }}y\right\vert
^{2}\leq \left\vert y\right\vert _{A_{R}}^{2}\leq \frac{4}{R\,\underline{\lambda} _{S}}%
\left\vert \Pi _{S}y\right\vert ^{2}+\frac{4}{R^{2}\,\underline{\lambda} _{S^{\bot }}}%
\left\vert \Pi _{S^{\bot }}y\right\vert ^{2}.  \label{Norm3a}
\end{equation}%
In particular, if $\left\vert A\right\vert _{\infty
}=\max_{i=1,...,m}\left\vert A_{i}\right\vert $ and $R\leq \underline{\lambda}
_{S}/4m\left\vert A\right\vert _{\infty }$ then
\begin{equation}
\frac{1}{4R\left\vert A\right\vert _{\infty }}\left\vert \Pi
_{S}y\right\vert ^{2}+\frac{1}{4R^{2}\left\vert A\right\vert _{\infty }}%
\left\vert \Pi _{S^{\bot }}y\right\vert ^{2}\leq \left\vert y\right\vert
_{A_{R}}^{2}\leq \frac{4}{R\,\underline{\lambda} _{S}}\left\vert \Pi _{S}y\right\vert ^{2}+%
\frac{4}{R^{2}\,\underline{\lambda} _{S^{\bot }}}\left\vert \Pi _{S^{\bot }}y\right\vert
^{2}.  \label{Norm3b}
\end{equation}
\end{proposition}

\textbf{Proof}. In a first stage we assume that $A_{i}\bot A_{j}$ for $i\leq
d<j.$ We will drop out this restriction in the second part of the proof. Let
$T_{S}$ and $T_{S^\bot }$ be the restriction of $y\mapsto
A_{R}A_{R}^{\ast }y$ to $S$ and to $S^{\bot }$ respectively. Since%
\begin{equation*}
A_{R}A_{R}^{\ast }y=R\sum_{i=1}^{d}\left\langle A_{i},y\right\rangle
A_{i}+R^{2}\sum_{i=d+1}^{m}\left\langle A_{i},y\right\rangle A_{i}
\end{equation*}%
our orthogonality hypothesis implies that $T_{S}y=R\sum_{i=1}^{d}\left%
\langle A_{i},y\right\rangle A_{i}\in S$ for $y\in S$ and $T_{S^\bot
}y=R^{2}\sum_{i=d+1}^{m}\left\langle A_{i},y\right\rangle A_{i}\in S^{\bot }$
for $y\in S^{\bot }.$ Since $A_{R}A_{R}^{\ast }$ is invertible it follows
that $T_{S}$ (respectively $T_{S^\bot })$ is an invertible operator from $S$
(respectively from $S^{\bot })$ into itself. For $y\in S$ we have
\begin{equation*}
R\,\overline{\lambda} _{S}\left\vert y\right\vert ^{2}\geq \left\langle
T_{S}y,y\right\rangle =R\sum_{i=1}^{d}\left\langle A_{i},y\right\rangle
^{2}\geq R\,\underline{\lambda} _{S}\left\vert y\right\vert ^{2}
\end{equation*}%
and since $\left\vert y\right\vert _{A_{R}}^{2}=\left\langle
T_{S}^{-1}y,y\right\rangle $ for $y\in S,$ we obtain
\begin{equation*}
\frac{1}{R\,\overline{\lambda} _{S}}\left\vert y\right\vert ^{2}\leq \left\vert
y\right\vert _{A_{R}}^{2}\leq \frac{1}{R\,\underline{\lambda} _{S}}\left\vert y\right\vert
^{2},\quad\mbox{if }y\in S.
\end{equation*}%
Similarly, we get
\begin{equation*}
\frac{1}{R^{2}\,\overline{\lambda} _{S^{\bot }}}\left\vert y\right\vert ^{2}\leq
\left\vert y\right\vert _{A_{R}}^{2}\leq \frac{1}{R^{2}\,\underline{\lambda} _{S^{\bot }}}%
\left\vert y\right\vert ^{2},\quad\mbox{if } y\in S^\bot.
\end{equation*}%
Let $y\in R^{n}.$ Since $(A_{R}A_{R}^{\ast })^{-1}\Pi _{S}y\in S$ we have $%
\left\langle (A_{R}A_{R}^{\ast })^{-1}\Pi _{S}y,\Pi _{S^{\bot
}}y\right\rangle =0$ so that $\left\vert y\right\vert
_{A_{R}}^{2}=\left\vert \Pi _{S}y\right\vert _{A_{R}}^{2}+\left\vert \Pi
_{S^{\bot }}y\right\vert _{A_{R}}^{2}.$ We obtain%
\begin{equation}
\frac{1}{R\,\overline{\lambda} _{S}}\left\vert \Pi _{S}y\right\vert ^{2}+\frac{1}{%
R^{2}\,\overline{\lambda} _{S^{\bot }}}\left\vert \Pi _{S^{\bot }}y\right\vert
^{2}\leq \left\vert y\right\vert _{A_{R}}^{2}\leq \frac{1}{R\,\underline{\lambda} _{S}}%
\left\vert \Pi _{S}y\right\vert ^{2}+\frac{1}{R^{2}\,\underline{\lambda} _{S^{\bot }}}%
\left\vert \Pi _{S^{\bot }}y\right\vert ^{2}.  \label{Norm5a}
\end{equation}

We drop now out the orthogonality assumption. For $j>d$ we consider the
decomposition $A_{j}=\Pi _{S}A_{j}+\Pi _{S^{\bot }}A_{j}$ and we define the
matrices $\overline{A}_{R}=(\sqrt{R}A_{1},...,\sqrt{R}A_{d},R\Pi _{S^{\bot
}}A_{d+1},...,$ $R\Pi _{S^{\bot }}A_{m})$ and $\widehat{A}_{R}=(0,...,0,R\Pi
_{S}A_{d+1},...,R\Pi _{S}A_{m})$ so that $A_{R}=\overline{A}_{R}+\widehat{A}%
_{R}.$ We will check that under the restriction (\ref{Norm2a}) we have%
\begin{equation}
4\left\vert \overline{A}_{R}^{\ast }y\right\vert ^{2}\geq \left\vert
A_{R}^{\ast }y\right\vert ^{2}\geq \frac{1}{4}\left\vert \overline{A}%
_{R}^{\ast }y\right\vert ^{2}\quad \forall y\in R^{n}.  \label{Norm4a}
\end{equation}%
We suppose for the moment that the above inequality is true and we prove (%
\ref{Norm3a}). Since $\left\vert A_{R}^{\ast }y\right\vert ^{2}=\left\langle
A_{R}A_{R}^{\ast }y,y\right\rangle $ the above inequality means that $4%
\overline{A}_{R}\overline{A}_{R}^{\ast }\geq A_{R}A_{R}^{\ast }\geq \frac{1}{%
4}\overline{A}_{R}\overline{A}_{R}^{\ast }$ and this gives
\begin{equation}
\frac{1}{4}\left\vert y\right\vert _{\overline{A}_{R}}^{2}\leq \left\vert
y\right\vert _{A_{R}}^{2}\leq 4\left\vert y\right\vert _{\overline{A}%
_{R}}^{2}.  \label{Norm6a}
\end{equation}%
Since the columns of $\overline{A}_{R}$ verify the orthogonality assumption
we may use the result from the first step with $A$ replaced with
$\overline{A}=(\overline{A}_{1},\ldots,\overline{A}_d,\overline{A}_{d+1},\ldots\overline{A}_m)$, with $\overline{A_j}=A_j$ if $j\leq d$ and $\overline{A_j}=\Pi_{S^{\bot }}A_{j}$ for $j>d$. Here, we have $\overline{S}=Vect\{\overline{A_1},\ldots,\overline{A_d}\}
=Vect\{A_1,\ldots,A_d\}=S$, so that $\underline{\lambda}_{\overline S}
=\underline{\lambda}_S$ and $\overline{\lambda}_{\overline S}
=\overline{\lambda}_S$. Moreover, since $\overline{S}^\bot=S^\bot$, the computations in (\ref{Norm1a})  are actually done with $%
\xi \in S^{\bot }$, and thus we obtain $\underline{\lambda }_{\overline{S}^{\bot }}=\underline{\lambda}
_{S^{\bot }}$ and $\overline{\lambda }_{\overline{S}^{\bot }}=\overline{\lambda} _{S^{\bot
}}.$ So (\ref{Norm5a}) gives%
\begin{equation*}
\frac{1}{R\,\overline{\lambda} _{S}}\left\vert \Pi _{S}y\right\vert ^{2}+\frac{1}{%
R^{2}\,\overline{\lambda} _{S^{\bot }}}\left\vert \Pi _{S^{\bot }}y\right\vert
^{2}\leq \left\vert y\right\vert _{\overline{A}_{R}}^{2}\leq \frac{1}{%
R\,\underline{\lambda} _{S}}\left\vert \Pi _{S}y\right\vert ^{2}+\frac{1}{R^{2}\,\underline{\lambda}
_{S^{\bot }}}\left\vert \Pi _{S^{\bot }}y\right\vert ^{2}
\end{equation*}%
which together with (\ref{Norm6a}) imply (\ref{Norm3a}).

It remains to prove (\ref{Norm4a}). We have%
\begin{equation*}
\left\vert \overline{A}_{R}^{\ast }y\right\vert ^{2}\geq
R\sum_{j=1}^{d}\left\langle A_{j},y\right\rangle
^{2}=R\sum_{j=1}^{d}\left\langle A_{j},\Pi _{S}y\right\rangle ^{2}\geq
R\,\underline{\lambda} _{S}\left\vert \Pi _{S}y\right\vert ^{2}
\end{equation*}%
and%
\begin{equation*}
\left\vert \widehat{A}_{R}^{\ast }y\right\vert
^{2}=R^{2}\sum_{j=d+1}^{m}\left\langle \Pi _{S}A_{j},y\right\rangle
^{2}=R^{2}\sum_{j=d+1}^{m}\left\langle \Pi _{S}A_{j},\Pi _{S}y\right\rangle
^{2}\leq R^{2}\left\vert \Pi _{S}y\right\vert ^{2}\sum_{j=d+1}^{m}\left\vert
\Pi _{S}A_{j}\right\vert ^{2}.
\end{equation*}%
Then (\ref{Norm2a}) gives $\left\vert \overline{A}_{R}^{\ast }y\right\vert
^{2}\geq 4\left\vert \widehat{A}_{R}^{\ast }y\right\vert ^{2}.$ Using the
inequality $(a+b)^{2}\geq \frac{1}{2}a^{2}-b^{2}$ we obtain%
\begin{equation*}
\left\vert A_{R}^{\ast }y\right\vert ^{2}=\left\vert \overline{A}_{R}^{\ast
}y+\widehat{A}_{R}^{\ast }y\right\vert ^{2}\geq \frac{1}{2}\left\vert
\overline{A}_{R}^{\ast }y\right\vert ^{2}-\left\vert \widehat{A}_{R}^{\ast
}y\right\vert ^{2}\geq \frac{1}{4}\left\vert \overline{A}_{R}^{\ast
}y\right\vert ^{2}
\end{equation*}%
and using $(a+b)^{2}\leq 2a^{2}+2b^{2}$ we get $\left\vert A_{R}^{\ast
}y\right\vert ^{2}\leq 4\left\vert \overline{A}_{R}^{\ast }y\right\vert
^{2}. $ $\square $

\smallskip

From now on we consider the specific situation when $a_{i}=\sigma
_{i}(t,x),[a]_{i,j}=[\sigma _{i},\sigma _{j}](t,x)$ and we denote by $A(t,x)$
respectively $A_{R}(t,x)$ the matrices associated to these coefficients. We
will need the following

\begin{lemma}\label{page-norm3}
\label{NORM3}
Let $x,y\in R^{n}$ be such that $|x-y|\leq 1$ and let $s,t\in [0,1]$. Assume that
\begin{equation}
\left\vert x-y\right\vert +\left\vert t-s\right\vert \leq \frac{\sqrt{%
\lambda _{\ast }(A(t,x))}}{(8dm)n^{2}(t,x)}\times \sqrt{\delta }
\label{Norm5}
\end{equation}%
Then for every $z\in R^{n}$ and $\delta\leq 1$ one has
\begin{equation}
\frac{1}{4}\left\vert z\right\vert^2 _{A_{\delta }(t,x)}\leq \left\vert
z\right\vert^2 _{A_{\delta }(s,y)}\leq 4\left\vert z\right\vert^2 _{A_{\delta
}(t,x)}.  \label{Norm6}
\end{equation}
\end{lemma}

\textbf{Proof}. The inequality (\ref{Norm6}) is equivalent to
\begin{equation*}
4(A_{\delta }A_{\delta }^{\ast })(t,x)\geq (A_{\delta }A_{\delta }^{\ast
})(s,y)\geq \frac{1}{4}(A_{\delta }A_{\delta }^{\ast })(t,x).
\end{equation*}%
We use the numerical inequality $(a+b)^{2}\geq \frac{1}{2}a^{2}-b^{2}$,
the hypothesis (\ref{Norm5}) and we obtain
\begin{eqnarray*}
\left\langle (A_{\delta }A_{\delta }^{\ast })(s,y)z,z\right\rangle
&=&\sum_{k=1}^{m}\left\langle A_{\delta ,k}(s,y),z\right\rangle
^{2}\\
&=&\sum_{k=1}^{m}(\left\langle A_{\delta ,k}(t,x),z\right\rangle
+\left\langle A_{\delta ,k}(s,y)-A_{\delta ,k}(t,x),z\right\rangle )^{2} \\
&\geq &\frac{1}{2}\sum_{k=1}^{m}\left\langle A_{\delta
,k}(t,x),z\right\rangle ^{2}-\sum_{k=1}^{m}(\left\langle A_{\delta
,k}(s,y)-A_{\delta ,k}(t,x),z\right\rangle )^{2} \\
&\geq &\frac{1}{2}\sum_{k=1}^{m}\left\langle A_{\delta
,k}(t,x),z\right\rangle ^{2}-(2dm)^{2}n^{4}(t,x)\delta (\left\vert
x-y\right\vert ^{2}+\left\vert t-s\right\vert ^{2})\times \left\vert
z\right\vert ^{2}.
\end{eqnarray*}%
Since $\lambda _{\ast }(A_{\delta }(t,x))\geq \delta ^{2}\lambda _{\ast
}(A(t,x))$ our hypothesis says that%
\begin{align*}
(2dm)^{2}n^{4}(t,x)\delta (\left\vert x-y\right\vert ^{2}+\left\vert
t-s\right\vert ^{2})\times \left\vert z\right\vert ^{2}
\leq &\frac{1}{4}%
\delta ^{2}\lambda _{\ast }(A(t,x))\times \left\vert z\right\vert ^{2} \\
\leq &\frac{1}{4}\lambda _{\ast }(A_{\delta }(t,x))\times \left\vert
z\right\vert ^{2}\\
\leq &\frac{1}{4}\sum_{k=1}^{m}\left\langle A_{\delta
,k}(t,x),z\right\rangle ^{2}
\end{align*}%
so that%
\begin{equation*}
\left\langle (A_{\delta }A_{\delta }^{\ast })(s,y)z,z\right\rangle \geq
\frac{1}{4}\sum_{k=1}^{m}\left\langle A_{\delta ,k}(t,x),z\right\rangle ^{2}=%
\frac{1}{4}\left\langle (A_{\delta }A_{\delta }^{\ast
})(t,x)z,z\right\rangle .
\end{equation*}%
Using $(a+b)^{2}\leq 2a^{2}+2b^{2}$ one proves the other inequality. $%
\square $

\medskip

In the last part of this section we establish the link between the norm $%
\left\vert z\right\vert _{A_{R}(t,x)}$ and the control (Caratheodory)
distance. We will use in a crucial way the alternative characterizations
given in \cite{[NSW]} for this distance - and these results hold in the homogeneous
case: the coefficients of the equations do not depend on time any more, so that we suppose now $\sigma
_{j}(t,x)=\sigma _{j}(x).$ Consequently, we handle the matrix $A_{R}(x)$
instead of $A_{R}(t,x).$

\smallskip

We first introduce a semi-distance $d$ on an open set $\Omega\subset R^n$ which is
naturally associated to the family of norms $\left\vert y\right\vert
_{A_{R}(x)}.$

\smallskip

We set $\Omega=\{x\in R^n\,:\,\lambda_*(A(x))>0\}=\{x\,:\,\det(AA^*(x))\neq 0\}$, which is open because $x\mapsto \det AA^*(x)$ is continuous. Notice that if $x\in\Omega$ then $\det A_RA^*_R(x)>0$ for every $R>0$.
For $x,y\in \Omega$, we define $d(x,y)$ by $d(x,y)<\sqrt{R}$ if and only if $%
\left\vert y-x\right\vert _{A_{R}(x)}<1.$ The motivation of taking $\sqrt{R}$
is the following: if we are in the elliptic case then $\left\vert
y-x\right\vert _{A_{R}(x)}\sim R^{-1/2}\left\vert y-x\right\vert $ so $%
\left\vert y-x\right\vert _{A_{R}(x)}\leq 1$ amounts to $\left\vert
y-x\right\vert \leq \sqrt{R}.$

It is straightforward to see that $d$ is a semi-distance on $\Omega$, in the sense that $d$ verifies the following three properties (see \cite{[NSW]}):\label{Norm7}%
\begin{itemize}
\item[$i)$]
for every $r>0$, the set $\{y\in\Omega\,:\,d(x,y)<r\}$ is open;
\item[$ii)$]
$d(x,y)=0$ if and only if $x=y$;
\item[$iii)$]
for every compact set $K\Subset\Omega$ there exists $C>0$ such that for every $x,y,z\in K$ one has $d(x,y)\leq C\big(d(x,z)+d(z,y)\big)$.
\end{itemize}
Moreover, one says that $d_{1}:\Omega\times \Omega\rightarrow R_{+}$ and $%
d_{2}:\Omega\times \Omega\rightarrow R_{+}$ are equivalent if for every
compact set $K\Subset\Omega$ there exists a constant $C$ such that for every
$x,y\in K$
\begin{equation}
\frac{1}{C}d_{1}(x,y)\leq d_{2}(x,y)\leq Cd_{1}(x,y).  \label{Norm7a}
\end{equation}%
In particular if $d_{1}$ is a distance and $d_{2}$ is equivalent with $d_{1}$
then $d_{2}$ is a semi-distance. And one says that $d_{1}$ is locally
equivalent with $d_{2}$ if for every $x_{0}\in \Omega$ there exists a
neighborhood $V$ of $x_{0}$ and a constant $C$ such that (\ref{Norm7a}) holds
for every $x,y\in V.$

We introduce now the control metric. For $x,y\in R^{n}$ we denote by $C(x,y)$
the set of controls $\psi \in L^{2}([0,1];R^{n})$ such that the
corresponding skeleton $du_{t}(\psi )=\sum_{j=1}^{d}\sigma _{j}(u_{t}(\psi
))\psi _{t}^{j}dt$ with $u_{0}(\psi )=x$ satisfies $u_{1}(\psi )=y.$ Notice
that the drift $b$ does not appear in the equation of $u_{t}(\psi ).$ Then
we define%
\begin{equation*}
d_{c}(x,y)=\inf \Big\{\Big(\int_{0}^{1}\left\vert \psi _{s}\right\vert
^{2}ds\Big)^{1/2}:\psi \in C(x,y)\Big\}.
\end{equation*}

\begin{theorem}
\label{NORM4}
\textbf{A}. Let
$$
\alpha (x)=\frac{\lambda _{\ast }^{1/2}(A(x))}{256d^{6}n^{6}(x)}.
$$
Then for every $x,y\in \Omega$ such that $d_{c}(x,y)\leq \frac{1}{4}\alpha ^{2}(x)$
one has $d(x,y)\leq 9\alpha ^{2}(x)d_{c}(x,y).$

\textbf{B}. $d$ is locally equivalent to $d_{c}$ on $\Omega$.

\textbf{C}. In particular for every compact set $K\Subset\Omega$ there exists $r_{K}$ and
$C_{K}$ such that for every $x,y\in K$ with $d(x,y)\leq r_{K}$ one has $%
d_{c}(x,y)\leq C_{K}d(x,y).$
\end{theorem}

\textbf{Proof A. }Let $\delta >0$, $x,y\in\Omega$ and $\psi \in C(x,y).$ Setting  $x_{t}=u_{t/\delta }(\psi )$, we obtain $dx_{t}=\sum_{j=1}^{d}%
\sigma _{j}(x_{t})\phi _{t}^{j}dt$ with $\phi
(t)=\delta ^{-1}\psi (t\delta ^{-1})$, which means that $x_{t}=u_{t}(\phi ).$
Notice also that $\int_{0}^{1}\left\vert \psi _{s}\right\vert ^{2}ds=\delta
\int_{0}^{\delta }\left\vert \phi _{s}\right\vert ^{2}ds.$

We denote now $C_{\delta }(x,y)$ the set of controls $\phi \in
L^{2}([0,\delta ];R^{n})$ such that the corresponding skeleton $u_{t}(\phi )$
with $u_{0}(\phi )=x$ verifies $u_{\delta }(\phi )=y.$ As a consequence
of the previous computations, one has
$$
d_{c}(x,y)=\sqrt\delta \times \inf \Big\{\Big(\int_{0}^{\delta }\left\vert \phi
_{s}\right\vert ^{2}ds\Big)^{1/2}:\phi \in C_{\delta }(x,y)\Big\}\equiv \sqrt{\delta }%
\times \inf \{\varepsilon _{\phi }(\delta ):\phi \in C_{\delta }(x,y)\}.
$$
Suppose now that $d_{c}(x,y)\leq \frac{1}{4}\alpha ^{2}(x).$ We take $\delta
=\frac{1}{4}\alpha ^{2}(x)$ so that $d_{c}(x,y)\leq \frac{1}{2}\alpha (x)%
\sqrt{\delta }.$ Then one may find a control $\phi \in C_{\delta }(x,y)$
such that $\varepsilon _{\phi }(\delta )\leq \frac{1}{2}\alpha (x)$ and $%
y=x_{\delta }(\phi ).$ Since $\varepsilon _{\phi }(\delta )+\sqrt{\delta }%
\leq \alpha (x)$ we may use  (\ref{Corection2}) and we obtain $\left\vert
y-x\right\vert _{A_{\delta }(x)}\leq 4\varepsilon _{\phi }(\delta )+\sqrt{%
\delta }\leq 3\alpha (x).$ It follows that%
\begin{equation*}
\left\vert y-x\right\vert _{A_{9\alpha ^{2}(x)\delta }(x)}\leq \frac{1}{%
3\alpha (x)}\left\vert y-x\right\vert _{A_{\delta }(x)}\leq 1
\end{equation*}%
and this gives $d(x,y)\leq 9\alpha ^{2}(x)\delta =9\alpha ^{2}(x)\times
\frac{1}{4}\alpha ^{2}(x).$ And this guarantees that $d(x,y)\leq 9\alpha
^{2}(x)\times d_{c}(x,y).$

\smallskip

\textbf{B.} We  prove now the converse inequality. We use the results
from \cite{[NSW]} so we recall the definition of the semi-distance $d_{\ast }$
(which is denoted by $\rho _{2}$ in \cite{[NSW]}). Given $\phi ^{i},\phi ^{k,j}\in
R,$ $i=0,...,d,1\leq k<j\leq d$ we consider the equation%
\begin{equation}
v_{t}(\phi )=x+\int_{0}^{t}(\sum_{j=1}^{d}\phi ^{j}\sigma _{j}(v_{s}(\phi
))+\sum_{i\neq j}\phi ^{i,j}[\sigma _{i},\sigma _{j}](v_{s}(\phi )))ds.
\label{Norm15}
\end{equation}%
Notice that $\phi ^{j},\phi ^{i,j}$ are now real constants (in contrast with
the time depending controls in the standard skeleton) and we have added the
vector fields $[\sigma _{i},\sigma _{j}]$ which does not appear in
skeletons. And the drift term $b$ does not appear. We denote by $P_{\ast
}(x,y)$ the family of paths $v_{t}(\phi )$ which satisfy (\ref{Norm15})
and such that $v_{1}(\phi )=y.$ We define $d_{\ast }$ by: $d_{\ast
}(x,y)<\delta $ if and only if one may find $\phi ^{i},\phi ^{k,j}\in R,$ $%
1\leq i,k,j\leq d,j\neq k$ such that $v_{\cdot}(\phi )\in P_{\ast }(x,y)$ and $%
\left\vert \phi ^{i}\right\vert <\delta ,\left\vert \phi ^{i,j}\right\vert
<\delta ^{2}.$ As a consequence of Theorem 2 and Theorem 4 from \cite{[NSW]} $%
d_{\ast }$ is locally equivalent with $d_{c}.$ So our aim is to prove that $%
d_{\ast }$ is locally dominated by $d.$ Let us be more precise: we fix $x\in
\Omega$ and we look for two constants $C_{x},\delta _{x}>0$ such that the
following holds: if $0<\delta \leq \delta _{x}$ and $d(x,y)\leq \sqrt{\delta
}$ then one may construct a control $\phi \in R^{m}$ such that $v_{\cdot}(\phi
)\in P_{\ast }(x,y)$ and $\left\vert \phi ^{i}\right\vert <C_{x}\sqrt\delta
,\left\vert \phi ^{i,j}\right\vert <C^2_{x}\delta .$ This implies $d_{\ast
}(x,y)\leq C_x\sqrt{\delta }$, and the statement will hold. Notice that we discuss local
equivalence, that is why we may take $C_{x},\delta _{x}$ depending on $%
x. $

We recall that $A_{i}(x),i=1,...,m$ is an enumeration of $\sigma
_{i}(x),[\sigma _{j},\sigma _{p}](x),i,j,p=1,...,d$ and that they span $R^{n}$ because $x\in\Omega$.
So, we choose $i_{1}<...<i_{d^{\prime }}\leq d<i_{d^{\prime
}+1}<...<i_{n}\leq m$ such that $A_{i_{k}}(x),k=1,...,d^{\prime }$ span $%
Vect\{A_{1}(x),...,A_{d}(x)\}$ and $A_{i_{k}}(x),k=1,...,n$ span $R^{n}.$ In
particular all of them are linearly independent. We denote $%
B_{k}(x)=A_{i_{k}}(x)$ and we want to use Theorem \ref{INVERSE3} for them.
Notice that $Vect\{A_{1}(x),...,A_{d}(x)\}=Vect\{B_{1}(x),...,B_{d^{\prime
}}(x)\}$ so the projections $\Pi $ and $\Pi ^{\bot }$ considered in Theorem %
\ref{INVERSE3} and in Proposition \ref{Norm2} coincide. In particular if $%
d(x,y)\leq \sqrt{\delta }$ then $\left\vert \Pi (y-x)\right\vert \leq
\left\vert A(x)\right\vert _{\infty }\sqrt{\delta }$ and $\left\vert \Pi
^{\bot }(y-x)\right\vert \leq \left\vert A(x)\right\vert _{\infty }\delta .$
And this also implies that $\left\vert y-x\right\vert \leq 2\left\vert
A(x)\right\vert _{\infty }\sqrt{\delta }.$

As $\theta\in R^n$, we look for a solution  to the equation
\begin{equation*}
y=\xi_1(\theta),\quad\mbox{with}\quad \xi_{t}(\theta )=x+\sum_{k=1}^{n}\theta _{k}\int_{0}^{t}B_{k}(\xi_{s}(\theta
))ds=x+\sum_{k=1}^{n}\theta _{k}\int_{0}^{t}A_{i_{k}}(\xi_{s}(\theta ))ds.
\end{equation*}%
So, we write it as%
\begin{equation*}
y=x+B(x)\theta +r(\theta )\quad \mbox{with}\quad r(\theta
)=\sum_{k=1}^{n}\theta _{k}\int_{0}^{t}(B_{k}(\xi_{s}(\theta ))-B_{k}(x))ds.
\end{equation*}%
Clearly $r\in C^{3}(R^{n},R^{n})$ and $r(0)=\nabla r(0)=0.$ Then,
\begin{equation*}
\left\vert y-x\right\vert \leq 2\left\vert A(x)\right\vert _{\infty }\sqrt{%
\delta }\leq 2\left\vert A(x)\right\vert _{\infty }\sqrt{\delta _{x}}
\end{equation*}%
and we suppose that $%
\delta _{x}$ is sufficiently small in order that $|y-x|$ satisfies (\ref{Comp3}), that is
$$
|y-x|<\frac {\lambda_*(B(x))^{1/2}}4\quad\mbox{and}\quad \left\vert y-x\right\vert < \frac{\lambda_*(B(x))}{8d^3(c_{2}(r)+c_{2}(r))}.
$$
We use then (\ref{Comp4}) and we obtain $\left\vert
\theta _{i}\right\vert \leq C_{x}\sqrt{\delta },$ $i=1,...,d^{\prime }$ and $%
\left\vert \theta _{i}\right\vert \leq C^2_{x}\delta ,$ $i=d^{\prime
}+1,...,n. $ This proves that $d_{\ast }(x,y)\leq C_{x}\sqrt{\delta }.$

\smallskip

\textbf{C.} For $x\in\Omega$, we denote $B_{d}(x,r):=\{y\in\Omega\,:\,d(x,y)<r\}$ and this is an open set.
Since $d$ and $d_{c}$ are locally equivalent for every compact $K\Subset\Omega$ and for every $x\in K$ there exists
$C_{x},\varepsilon _{x}>0$ such that for $y\in B_{d}(x,\varepsilon _{x})$ we
have $d_{c}(x,y)\leq C_{x}d(x,y).$ Since the set $K$ is compact we may find $%
x_{1},...,x_{N}\in K$ such that $K\subset \cup _{i=1}^{N}B_{d}(x_{i},\varepsilon
_{x_{i}}).$ We denote $C_{\max }=\max_{i=1,...,N}C_{x_{i}}.$ Let us prove
that there exists $r_{\ast }>0$ such that for every $x\in K$ and every $y\in
B_{d}(x,r_{\ast })$ we have $d_{c}(x,y)\leq C_{\max }d(x,y).$

For $x\in K$ one may find $i$ such that $x\in B_{d}(x_{i},\varepsilon
_{x_{i}})$ and $r>0$ such that $B_{d}(x,r)\subset B_{d}(x_{i},\varepsilon
_{x_{i}}).$ We define $r_{x}=\sup \{r>0:$ $\exists $ $i\in \{1,...,N\}$ such
that $B_{d}(x,r)\subset B_{d}(x_{i},\varepsilon _{x_{i}})\}.$ We claim that $%
r_{\ast }:=\inf_{x\in K}r_{x}>0.$ Indeed suppose that this is not true. Then
one may find a sequence $y_{n}\rightarrow y_{0}$ such that $%
r_{y_{n}}\rightarrow 0.$ Since $r_{y_{0}}>0$ there exists $n_{\ast }$ such
that for $n\geq n_{\ast }$ one has $B_{d}(y_{n},\frac{1}{2C_{K}}%
r_{y_{0}})\subset B_{d}(y_{0},r_{y_{0}})\subset B_{d}(x_{i},\varepsilon
_{t_{i}})$ for some $i.$ Here $C_{K}$ is the constant in the triangle
inequality $iii)$ at page \pageref{Norm7}. And this means that $r_{y_{n}}\geq \frac{1}{%
2}r_{y_{0}}>0$ which is in contradiction with our hypothesis. So we have
proved that $r_{\ast }>0.$

Consider now $y\in B_{d}(x,r_{\ast }).$ There exists $i$ such that $%
B_{d}(x,r_{\ast })\subset B_{d}(x_{i},\varepsilon _{x_{i}})$ and this means
that $y,x\in B_{d}(x_{i},\varepsilon _{x_{i}})$ and consequently $%
d_{c}(x,y)\leq C_{x_{i}}d(x,y)\leq C_{\max }d(x,y).$ $\square $

\medskip

Finally we give:

\medskip

\textbf{Proof of Proposition} \ref{NOT2}. We will first prove that under our
hypothesis $d(x,y)\leq \sqrt{\lambda _{x}/(4m)n^{4}(x)}.$ Let $R$ be such
that $d(x,y)\geq \sqrt{R}$ so that $\left\vert y-x\right\vert
_{A_{R}(x)}\geq 1.$ Then by (\ref{Norm3})
\begin{equation*}
\frac{\lambda _{x}\sqrt{\lambda _{\ast }(A)}}{(4m)n^{4}(x)}\geq \left\vert
y-x\right\vert \geq R\sqrt{\lambda _{\ast }(A)}\left\vert y-x\right\vert
_{A_{R}}\geq R\sqrt{\lambda _{\ast }(A)}.
\end{equation*}%
It follows that $R\leq \lambda _{x}/(4m)n^{4}(x)$ which proves our assertion.

We suppose now that $d(x,y)\geq \sqrt{R}.$ Since $R\leq \lambda
_{x}/(4m)n^{4}(x)$ we may use (\ref{Not10}) and we obtain
\begin{equation*}
\frac{4}{R\lambda _{x}}\left\vert \Pi _{x}(y-x)\right\vert ^{2}+\frac{4}{%
R^{2}\lambda _{x}^{\bot }}\left\vert \Pi _{x}^{\bot }(y-x)\right\vert
^{2}\geq \left\vert y-x\right\vert _{A_{R}(x)}\geq 1
\end{equation*}%
which gives $\overline{d}(x,y)\geq \sqrt{R}.$ So $\overline{d}(x,y)\geq
d(x,y).$ Suppose that $d(x,y)<\sqrt{R}.$ Then%
\begin{equation*}
\frac{1}{4Rn^{2}(x)}\left\vert \Pi _{x}(y-x)\right\vert ^{2}+\frac{1}{%
4R^{2}n^{2}(x)}\left\vert \Pi _{x}^{\bot }(y-x)\right\vert ^{2}\leq
\left\vert y-x\right\vert _{A_{R}(x)}^{2}<1
\end{equation*}%
and this reads $\underline{d}(x,y)<\sqrt{R}$ which gives $\underline{d}%
(x,y)\leq d(x,y).\square $

\end{document}